\documentclass[12pt,a4paper]{article}

\usepackage[english]{babel}
\usepackage[T1]{fontenc}
\usepackage{amssymb}
\usepackage{enumitem}
\usepackage{amsthm}
\usepackage{amsmath}
\usepackage[]{algorithm2e}
\usepackage{thmtools}
\usepackage{lscape}
\usepackage{cite}

\theoremstyle{definition}
\newtheorem{definition}{Definition}
\theoremstyle{remark}

\renewcommand\thmcontinues[2]{Continued}

\newcommand{\glebdone}[1]{}

\newcommand{\wudone}[1]{}

\usepackage[a4paper,top=3cm,bottom=3cm,left=2.5cm,right=2.5cm,marginparwidth=1.75cm]{geometry}

\usepackage{amsmath}
\usepackage{graphicx}
\usepackage{subfigure}
\usepackage[colorinlistoftodos]{todonotes}
\usepackage[colorlinks=true, allcolors=blue]{hyperref}

\title{Some ideas about graphic representations of discrete fuzzy measures}
\usepackage{authblk}
\author[1*]{Jian-Zhang Wu}
\author[1]{Gleb Beliakov}
\affil[1]{School of Information Technology, Deakin University, Burwood 3125, Australia}
\affil[*]
{Corresponding author: Jian-Zhang Wu, sjzwjz@gmail.com}

\begin{document}

\maketitle

\begin{abstract}
Graphs serve as efficient tools for visualizing mathematical concepts and their interrelationships. 
In this paper, focusing on the discrete case with universal set with finite elements, we first introduce the rules and characteristics of graph representation of fuzzy measure and discuss graphic properties of fuzzy measure's duality, symmetry, nonadditivity and nonmodularity. Then we show the graphic presentations of some special families of fuzzy measures, such as the k-additive measure, the k-maxitive and minitive measure, k-order representative measure, k-order interactive measure and the p-symmetric fuzzy measure, as well as of three nonlinear integrals, i.e., the Choquet integral, the Sugeno integral and the pan integral. 
Finally, we provide visualizations for the fuzzy measure fitting procedure and tools for comparing fuzzy measures.

\noindent \textbf{keywords}: Graph; Fuzzy measure; Nonlinear integral; fuzzy measure identification; Decision making.
\end{abstract}

\section{Introduction}

The discrete fuzzy measure, also known as a capacity, as introduced by Sugeno \cite{sugeno1974theory} and Choquet \cite{choquet1954}, is a set function defined on finite sets that exhibits normality and monotonicity properties.
Research in the field of fuzzy measures encompasses many branches, including investigations into equivalent linear transformations of fuzzy measures \cite{beliakov2019discrete_book,grabisch2016:setfunctionbook}, exploration of specific families of fuzzy measures, the development of relevant nonlinear integrals, and methods for identifying or fitting fuzzy measures.
In recent work, we have developed matrix representations \cite{wu2019matrixrepresentation} and marginal representations \cite{wu2019marginal} for fuzzy measure theories. This paper aims to introduce a graphical representation approach to further advance the field.

Through the graphical representation of fuzzy measures, where subsets serve as vertices with fuzzy measure values depicted as heights or radii of points and inclusion relationships represented as edges, we can elucidate the following key aspects:
\footnotemark

\footnotetext{Certainly, these ideas in this paper primarily occurred and finished between 2019 and 2020. Despite our efforts to refine and enhance them, we have encountered challenges in achieving the desired level of rigor. Now, we have made the decision to share these ideas with the public, with the hope that readers and peers can offer valuable insights and contribute to the ongoing refinement and enhancement of these initial thoughts.}
\begin{itemize}
    \item \textbf{Monotonicity:} The property of monotonicity in fuzzy measures is readily discernible as we observe the vertices of subsets exhibiting non-decreasing heights or radii, along with all edges having nonnegative heights.

 \item \textbf{Duality, Additivity, and Symmetry:} Concepts such as duality, additivity, and symmetry within fuzzy measures manifest as specific symmetric structures within the graphical representation.

 \item \textbf{Nonadditivity and Nonmodularity}: Nonadditivity and nonmodularity can be comprehended as dynamic trends in the heights of vertices associated with differing cardinalities.

 \item \textbf{Linear Equivalent Transformations:} Linear equivalent transformations of fuzzy measures result in flexible adjustments of the graphical structures, offering insights into their behavior.

 \item \textbf{Particular Families of Fuzzy Measures:} Various families of fuzzy measures tend to accentuate or amplify lower-level structural characteristics, often leading to distinctive layouts at higher levels within the graphical representation.

 \item \textbf{Nonlinear Integrals:} The mechanisms of aggregation in nonlinear integrals become more comprehensible when we examine their fundamental functions within the fuzzy measure graphs. The aggregation values and trends become apparent through the analysis of randomly generated fuzzy measures.

 \item \textbf{Fuzzy Measure Fitting:} The process of fitting a fuzzy measure becomes both visually and analytically accessible to decision-makers as they observe changes and evolutions in the graphical representation, aiding in the decision-making process.

 \item \textbf{Comparison about fuzzy measures:} 
 The analysis of subset similarity and clustering can be conducted using various indices, including M"obius transformation, nonadditivity index, and nonmodularity index. Furthermore, a comprehensive assessment of fuzzy measure characteristics can be achieved through indices such as entropy and orness.

\end{itemize}

The rest of this paper is organized as follows. We introduce the graphic representation of fuzzy measure and explain some basic properties of its linear transformations in Section 2. In Section 3 we discuss some particular families of fuzzy measures. Section 4 and 5 are for the graphic explanations of nonlinear integral and the fuzzy measure identification methods respectively. Section 6 is for the comparison analyses about subsets and whole fuzzy measure. Finally, we conclude the paper in Section 7.


\section{Fuzzy measures and the graphic representations of their basic characteristics}

Let $N = \{ 1,2, \ldots ,n\}$, $n \ge 2$, be the universal set of finite elements, $\mathcal{P}(N)$ be its power set, and $|S|$ be the cardinality of a subset $S\subseteq N$.
\theoremstyle{definition}
\begin{definition} \label{def. fuzzy measure} \cite{choquet1954,grabisch2008review,sugeno1974theory}
 A fuzzy measure on $N$ is a set function $\mu :\mathcal {P}(N) \to [0,1]$
such that
\begin{enumerate}
\item[(i)]{$\mu (\emptyset ) = 0, \mu (N) = 1;$ (boundary condition)}
\item[(ii)]{$\forall A,B \subseteq N$, $A \subseteq B$ implies $\mu (A) \le \mu (B)$. (monotonicity condition)}
\end{enumerate}
\end{definition}
The boundary conditions aim to normalize measure values, while monotonicity ensures that adding new elements strictly does not decrease the original coalition's measure value.

\begin{figure}[!htp]
	\color{black}
	\centering	
	\subfigure[The topological graph of $\mu_1$.]{
		\begin{minipage}[t]{0.5\linewidth}
			\centering
			\includegraphics[width=3in]{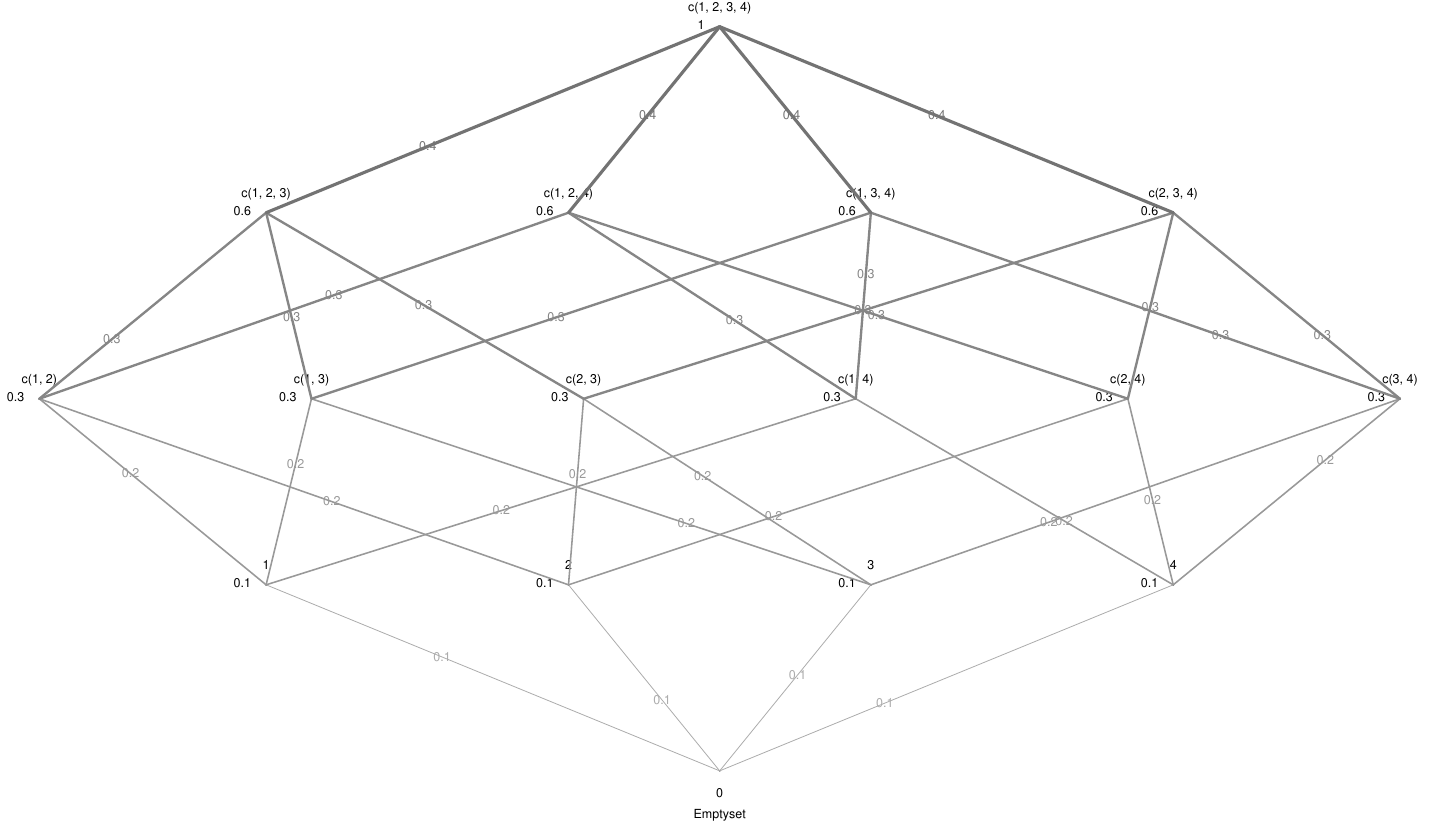}
		\end{minipage}%
	}%
	\subfigure[The "height-on" graph of $\mu_1$.]{
		\begin{minipage}[t]{0.5\linewidth}
			\centering
			\includegraphics[width=3in]{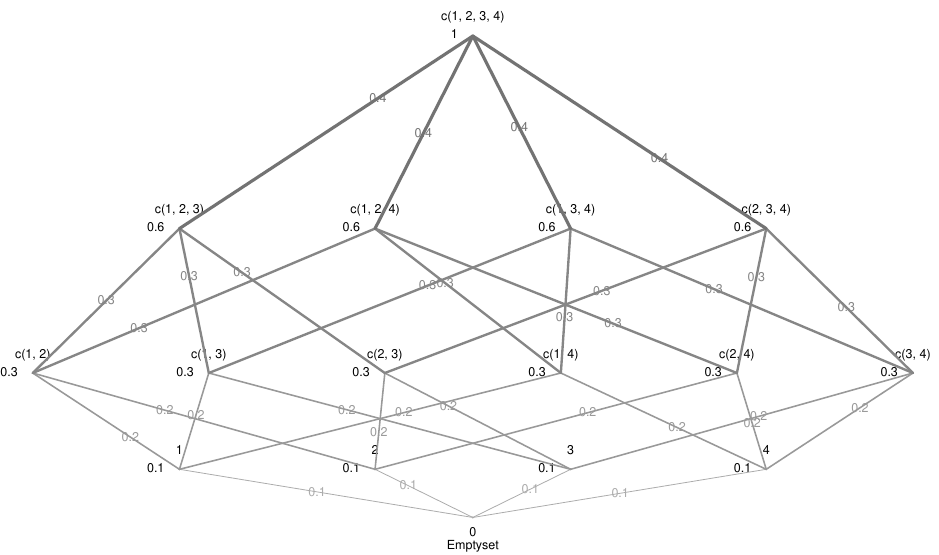}
		\end{minipage}%
	}
	\caption{Two types of graphic representations of fuzzy measure $\mu_1$ on $\{1,2,3,4\}$.} 
	\label{fig-4-symmetric-capa}
\end{figure}

A fuzzy measure can be visually represented as a graph in which subsets serve as points, inclusion relationships among them are depicted as edges, and fuzzy measure values are indicated by the heights of the points. For instance, in Figure \ref{fig-4-symmetric-capa}, we present the graph of a fuzzy measure, denoted as $\mu_1$, defined on the set $\{1,2,3,4\}$, where "c(1, 2)" denotes the coalition of criteria 1 and 2, and the value 0.45 represents its fuzzy measure value, all subsets are arranged in a specific topological order.

Here, we provide detailed rules governing the structure of the graph:
\begin{itemize}
	\item \textbf{Height:} we assign the height of the empty set point as 0 and the height of the universal criteria set point as 1, as specified in the standard definition of a fuzzy measure (see Definition \ref{def. fuzzy measure}). In the topological graph of a fuzzy measure, such as the one depicted in (a) of Fig. \ref{fig-4-symmetric-capa}, the height of each subset is determined by its cardinality. For instance, the height of subset $A$ is set to $|A|/n$. However, in certain situations, we opt to set the height of a subset to match its fuzzy measure value, facilitating a more intuitive understanding of the fuzzy measure. This type of graph is referred to as the "height-on" type, as shown in (b) of Fig. \ref{fig-4-symmetric-capa}. In this representation, you can observe that the larger the subset, the greater the increase in its height compared to its lower neighbors. The fuzzy measure value or height of a subset is positioned to the left of its corresponding point.
	
	\item \textbf{Horizontal Position:} in terms of horizontal positioning, the empty set and the universal set are always centered within the figure. Other subsets are evenly distributed and centered along a line that passes through the empty set and the universal set. Specifically, for a set of $n$ criteria, we assign the subsets with cardinality $\text{round}(n/2)$ the largest horizontal length $l$. Subsets with cardinality $r, r=1,...,n$  are allocated a horizontal length of ${\binom{n}{r}l}/{\binom{n}{\text{round}(n/2)}}$, where $\text{round}(x)$ returns the nearest integer to $x$.
	\item \textbf{Order:} we suggest adopting a dimension-based order, as seen in Figure \ref{fig-4-symmetric-capa}, where the point representing a subset $A$ is symmetric with its complement subset $N \setminus A$ in the horizontal perspective.
	\item \textbf{Monotonicity:} since each edge represents an inclusion relationship between two subsets, monotonicity becomes evident by observing the heights of subsets along chains within the graph. We can examine the maximum chains that start with the empty set and end with the universal set. For instance, consider the chain: "Empty set $\rightarrow$ 1 $\rightarrow$ c(1, 2) $\rightarrow$ c(1, 2, 3) $\rightarrow$ c(1, 2, 3, 4)." This chain reflects the inclusion relationships: $\emptyset \subset \{1\} \subset \{1,2\} \subset \{1,2,3\} \subset \{1,2,3,4\}$. The gradual increase in the heights of these subsets signifies that $ \mu(\emptyset) \leq  \mu(\{1\}) \leq  \mu(\{1,2\}) \leq  \mu(\{1,2,3\}) \leq  \mu(\{1,2,3,4\})$.
	\item \textbf{Marginal contribution:} 
 each edge can be considered as a marginal contribution, like edge "1 $\rightarrow$ c(1, 2)" can be seen as $\Delta_2\mu(\{1\})$, whose value is just the height of its edge, i.e., $\mu(\{1,2\})-\mu(\{1\})$. The marginal contribution is reflected by the width of its edge as well as its precise value on the edge, and also with a gradient color from "gray" to "black" (gray reflects 0 and black 1).

\end{itemize} 

Some characteristics of the fuzzy measure graph with $n$ criteria are give as follows:
\begin{itemize}
	\item Each point has $n$ edges. More specifically, the point representing set $A$ has $\binom{|A|}{|A|-1}$=$|A|$ edges connecting to neighbor subsets in the lower level, if they exist, and $n-|A|$ edges connecting to neighbor subsets in the upper level, if they exist.
	\item There are $n!$ maximum chains, and each maximum chain consists of $n+1$ points. The sum of all the vertical heights of edges in a maximum chain is always 1, which adheres to the boundary condition of a fuzzy measure where the fuzzy measure of the universal set is set to 1.
	\item The monotonicity condition can be just equivalently represented as all the $n!$ edges' heights are nonnegative. That is, monotonicity means the nondecreasing of the subsets' heights in each maximum chain according to their cardinalities.
\end{itemize}

\begin{definition} \label{def-dual-capa} \cite{beliakov2019discrete_book}
	The dual of fuzzy measure $\mu $ is defined as
	$$\overline{\mu}(A)=\mu(N)-\mu(N\setminus A).$$
\end{definition}

\begin{figure}[!htp]
	\color{black}
	\centering
	\includegraphics[width=0.6\textwidth]{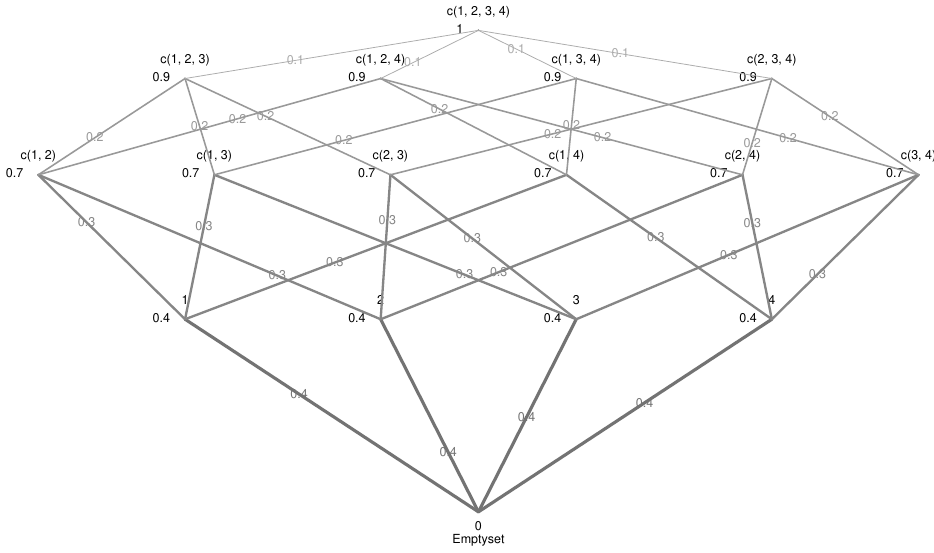}
	\caption{The "height-on" graph of $\overline{\mu}_1$.} 
	\label{fig-4-symmetric-capa-dual}
\end{figure}

The dual of the fuzzy measure $\mu_1$, denoted as $\overline{\mu}_1$ and shown in Figure \ref{fig-4-symmetric-capa}, is presented in Figure \ref{fig-4-symmetric-capa-dual}. One can easily observe that in both figures, the sum of the height of subset $A$ in one figure and the height of its complement subset $N \setminus A$ in the other figure is always 1. If a symmetric order, similar to that in these two figures, is adopted, we can further notice that the point representing a subset in one figure is symmetric with its complement subset in the other figure with respect to the middle point of the line segment between the points of the universal set and the empty set. For example, observe the point representing subset $\{1\}$ in Figure \ref{fig-4-symmetric-capa} and the point representing subset $\{2,3,4\}$ in Figure \ref{fig-4-symmetric-capa-dual}.


\begin{definition} \cite{beliakov2019discrete_book}
	A fuzzy measure $\mu$ on $N$ is said to be symmetric if its fuzzy measure values are based on the cadinality, i.e., if $|A|=|B|$ then $\mu(A)=\mu(B)$.
\end{definition}
Actually the Figs. \ref{fig-4-symmetric-capa}  and \ref{fig-4-symmetric-capa-dual} show two symmetric fuzzy measures, their points of subsets are really symmetric w.r.t. the line segment from the point of universal set to the point of empty set as well as its center point.

\begin{definition} \cite{beliakov2019discrete_book}
	A fuzzy measure $\mu$ on $N$ is said to be additive if $\mu(A \cup B)=\mu(A)+\mu(B)$ for any $A, B \subseteq N$ and $A \cap B = \emptyset$. Furthermore, $\mu$ is superadditive (subadditive)  if the preceding equation's direction is $\geq$ ($\leq$).
\end{definition}

\begin{figure}[!htp]
	\color{black}
	\centering
	\includegraphics[width=0.6\textwidth]{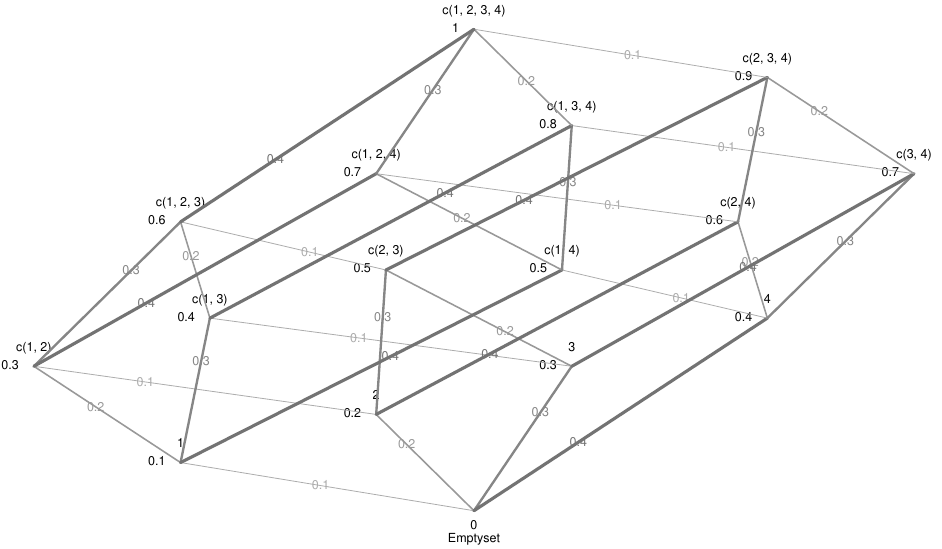}
	\caption{An "height-on" graph of additive measure $\mu_2$.} 
	\label{fig-4-additive-capa}
\end{figure}

It can be proved that the additive measure is self dual, hence in its graph, the subset is symmetric with its complement subset w.r.t. the center point of the line segment from the point of universal set to the point of empty set, even each edge has its symmetric edge w.r.t. the center point, see Figure \ref{fig-4-additive-capa}. 

\color{black}
We further introduce a concept called comprehensive importance or contribution of subset, which can be seen as an extension of concept of value \cite{shapley1953value}, comprehensive importance or contribution of singleton. 

\begin{definition}
	\label{def-shapley-importance}
	The Shapley comprehensive importance or contribution index of a subset $A \subseteq N$ w.r.t. $\mu $ is defined as
	$$k_\mu (A) = \sum\limits_{B \subseteq N\backslash A} \frac{1}{{|N| - |A| + 1}}{{\left( {\begin{array}{*{20}{c}}
	   {|N| - |A|}  \\
	   {|B|}  \\
	\end{array}} \right)}^{ - 1}}[\mu(A\cup B)-\mu (B)].$$
\end{definition}
One can figure out that $k_\mu(\{i\})$ is just the Shapley value of criterion $i$, hence comprehensive importance index is an extension of value. We always have $k_\mu(\emptyset)=0$,  $k_\mu(N)=1$ and $0 \leq k_\mu(A) \leq 1, \forall A \subseteq N.$ 

It can be verified that for additive measure $\mu$, 
\begin{equation*}
k_\mu(A)=\mu(A)=\sum_{i \in A} k_\mu(\{i\})=\sum_{i \in A} \mu(\{i\}), 
\end{equation*}
\begin{equation*}
k_\mu(N)=\mu(N)=\sum_{i \in N} k_\mu(\{i\})=\sum_{i \in N} \mu(\{i\})=1.
\end{equation*}

\begin{figure}[!htp]
	\color{black}
	\centering
	\includegraphics[width=0.6\textwidth]{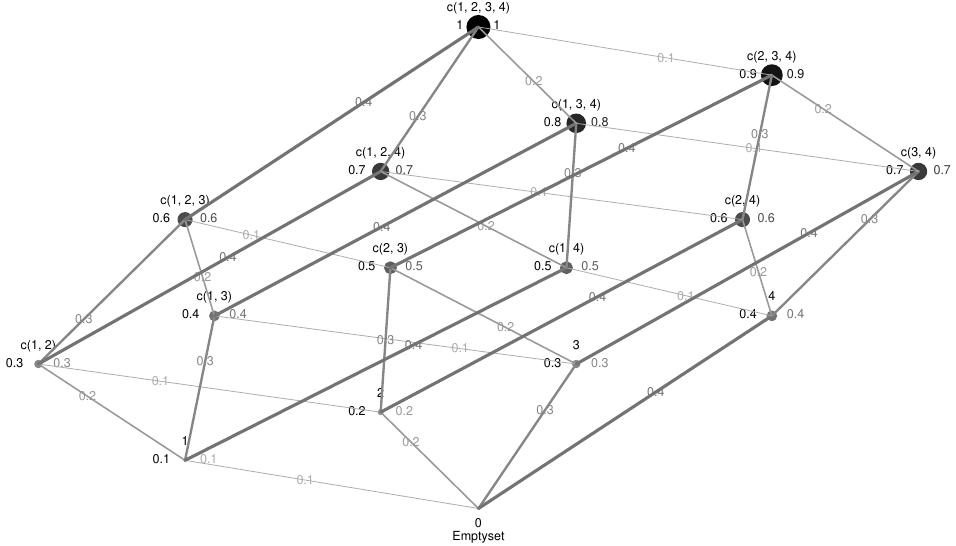}
	\caption{The "height-on" graph of $\mu_2$ with comprehensive importance.} 
	\label{fig-4-additive-capa-with-sh-importance}
\end{figure}

We can further represent the comprehensive importance of subsets using circles, where the radius is proportional to their importance values and displayed in a gradient color scale from gray to black. For example, in Fig. \ref{fig-4-additive-capa-with-sh-importance}, the fuzzy measure, as shown in Fig. \ref{fig-4-additive-capa}, is depicted, allowing one to easily verify the two equations associated with this additive measure.

Actually, it's easy to prove that Shapley comprehensive importance or contribution index of a fuzzy measure is still a fuzzy measure, and Shapley comprehensive importance or contribution index of an additive measure is still an additive measure.

\color{black}
\begin{definition} \label{def. non modular fuzzy measure }
	\cite{beliakov2019discrete_book,pap1995null}
	A fuzzy measure $\mu $ on $N$ is said to be supermodular (or submodular), if $\mu (A \cup B) + \mu (A \cap B) \geq (\text{or} \leq)  \mu (A) + \mu (B),$  $\forall A,B \subseteq N$.
\end{definition}
A fuzzy measure $\mu $ on $N$ is supermodular (submodular) if and only if 
\begin{equation*}
\Delta_i \mu(A) \geq ( \leq) \Delta_i \mu(B), \forall B \subset A, A, B \subseteq N, 
\end{equation*}
where $\Delta_i \mu(A)= \mu(A\cup \{i\})-\mu(A)$. 

\color{black}

According to this rule, we can observe that Figure \ref{fig-4-symmetric-capa} represents a supermodular fuzzy measure. For instance, consider the nondecreasing heights of edges: $\Delta_1 \mu(\emptyset)$,  $\Delta_1 \mu(\{2\})$,  $\Delta_1 \mu(\{2,3\})$,  $\Delta_1 \mu(\{2,3,4\})$. It's worth noting that these edges, or marginal contributions, are associated with the maximum chain "empty set $\rightarrow$ 2 $\rightarrow$ c(2, 3) $\rightarrow$ c(2,3,4) $\rightarrow$ c(1,2,3,4)" and criterion 1. This approach provides a maximum chain-based method for checking the nonmodularity of a fuzzy measure, especially through its graphical representation. Using this approach, it becomes clear that Figure \ref{fig-4-symmetric-capa-dual} represents a submodular fuzzy measure (the marginal contributions of singletons are non-ascending, or more precisely, decreasing), and Figure \ref{fig-4-additive-capa-with-sh-importance} represents a modular/additive measure (all the marginal contributions of a given criterion are the same).

%



These rigorous but time consuming ways to check the nonadditivity and nonmodularity is not efficient and almost impossible with a relatively large number of criteria. Hence, we introduce two types of indices to more easily check nonadditivity and nonmodularity, which is accepted from the perspective of probabilistic expectation in practice. 


\begin{definition} \label{def. some representation of fuzzy measure} \cite{grabisch2008review,wuBeliakovNonadd,wuBeliakovNonmodu}
Let $\mu $ be a fuzzy measure on $N$, 
the nonadditivity index of subset $A$ w.r.t. $\mu$ is defined as
\begin{equation*} \label{eq. nonaddi of fuzzy measure}
{n_\mu }(A) = \mu (A) - \frac{1}{{{2^{|A| - 1}} - 1}} \sum\limits_{C \subset A} {\mu (C)},
\end{equation*}
the nonmodularity index of subset $A \subseteq N$ w.r.t. $\mu$ is defined as
\begin{equation*} \label{eq-explicit-nonmodularity-index}
{d_\mu }(A) =\mu (A) -\frac{1}{|A|}\sum_{\{i\} \subset A } [\mu(\{i\})+\mu(A\backslash\{i\})].
\end{equation*}
\end{definition}

These two types of indices possess favorable properties for characterizing the interactions among decision criteria, such as a unified range, extreme interaction property, and moderate dummy property. Additional details can be found in \cite{wuBeliakovNonadd, wuBeliakovNonmodu, wuBeliakovproba-bi-}.

\begin{figure}[!htp]
	\color{black}
	\centering	
	\subfigure[$\mu_1$ with nonadditivity index.]{
		\begin{minipage}[t]{0.5\linewidth}
			\centering
			\includegraphics[width=3in]{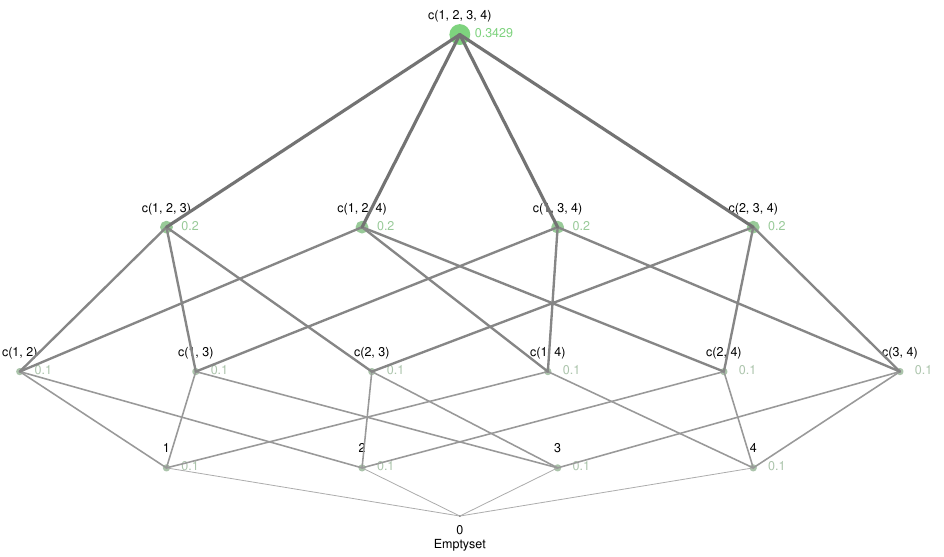}
		\end{minipage}%
	}%
	\subfigure[${\mu}_1$ with nonmodularity index.]{
		\begin{minipage}[t]{0.5\linewidth}
			\centering
			\includegraphics[width=3in]{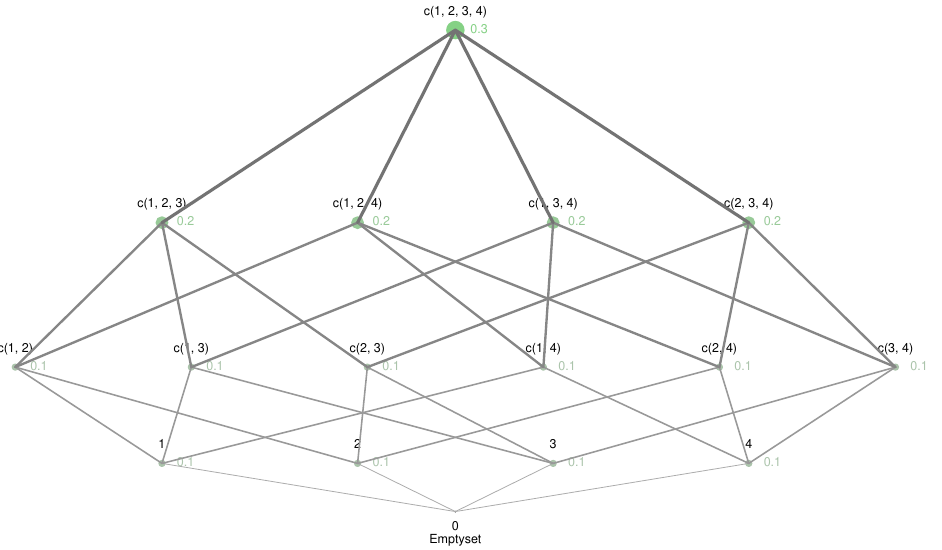}
		\end{minipage}%
	}
	
	\subfigure[$\overline{\mu}_1$ with nonadditivity index.]{
		\begin{minipage}[t]{0.5\linewidth}
			\centering
			\includegraphics[width=3in]{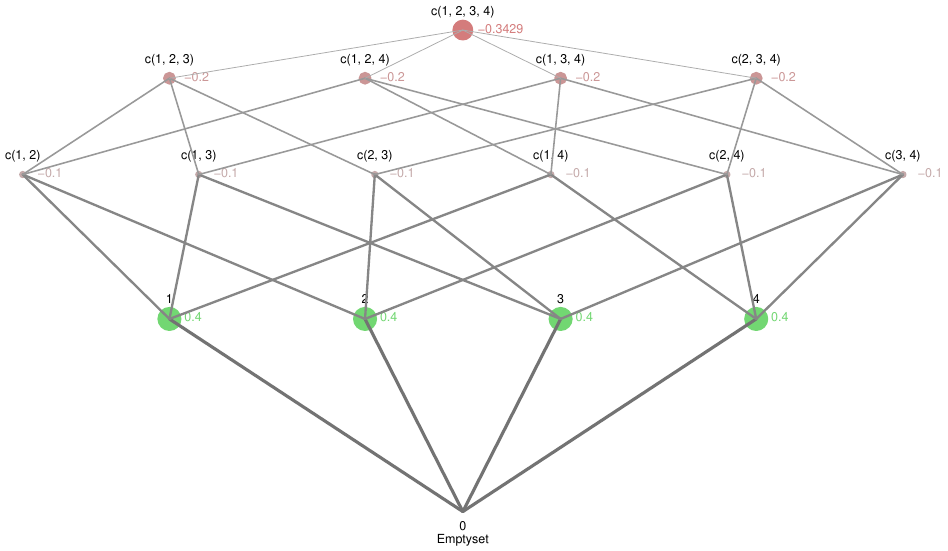}
		\end{minipage}%
	}%
	\subfigure[$\overline{\mu}_1$ with nonmodularity index.]{
		\begin{minipage}[t]{0.5\linewidth}
			\centering
			\includegraphics[width=3in]{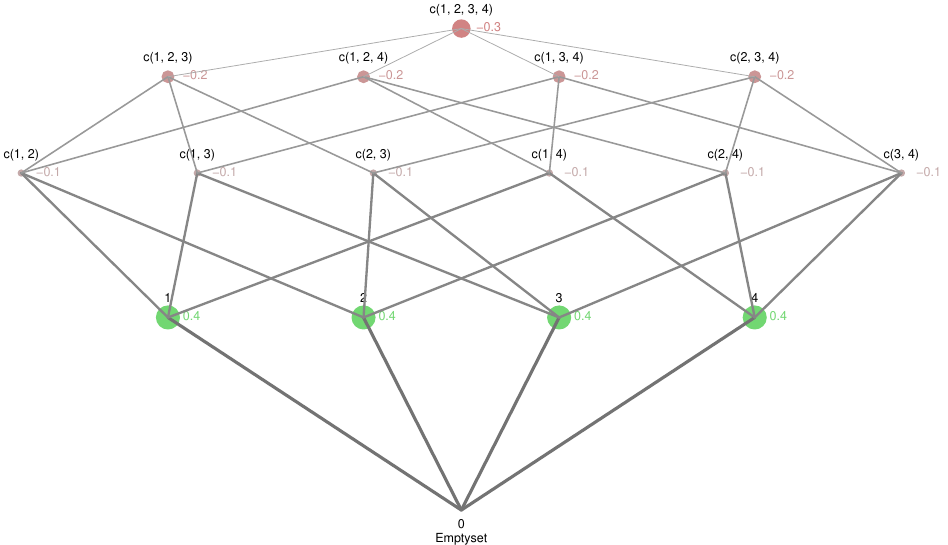}
		\end{minipage}%
	}
%
	\caption{The "height-on" graphs with nonadditivity and nonmodularity indices.}
	\label{fig-4-symmetric-and-dual-non-add-modu}
\end{figure}

Figure \ref{fig-4-symmetric-and-dual-non-add-modu} illustrates the nonadditivity and nonmodularity indices of fuzzy measures $\mu_1$ and $\overline{\mu}_1$. In the figure, points are sized based on the absolute values of the index values and are filled with a gradient color ranging from red to gray to green, where these three colors correspond to -1, 0, and 1, respectively. Additionally, the index values are displayed on the left side in corresponding colors.

From subfigures (a) and (b) in Fig. \ref{fig-4-symmetric-and-dual-non-add-modu}, it is evident that the nonadditivity and nonmodularity indices of nonsingleton and nonempty sets are all positive. 
Conversely, for subfigures (c) and (d), those values are all the negative.

 \color{black}


\section{Particular fuzzy measures and their graphs}

The primary objective of specific families of fuzzy measures is to effectively address the exponential complexity inherent in fuzzy measures on $n$ criteria, which involve $2^n-2$ coefficients. These families include the $k$-additive measure \cite{mesiar1999k}, the $k$-tolerant and intolerant fuzzy measure \cite{marichal2007k}, the $k$-maxitive and minitive measure \cite{mesiar2018k, wuBeliakovkminitive}, and the $p$-symmetric fuzzy measure \cite{miranda2002p}.

\begin{definition} \cite{grabisch1997k}
Let $k \in \{1,..., n\}$, a fuzzy measure $\mu$ on $N$ is said to be $k$-additive if its Möbius representation satisfies $m_\mu(A) = 0$ for all $A \subseteq N$ such that $|A| > k$ and there exists at least one subset $A$ of $k$ elements such that $m_\mu(A) \neq 0$, where the Möbius transformation of $\mu $ is defined as
$${m_\mu }(A) = \sum_{C \subseteq A} {{{( - 1)}^{|A\backslash C|}}\mu (C)}.$$
\end{definition}

\begin{figure}[!htp]
	\color{black}
	\centering	
	\subfigure[2-additive measure $\mu_3$ on $\{1,2,3,4\}$]{
		\begin{minipage}[t]{0.4\linewidth}
			\centering
			\includegraphics[width=3in]{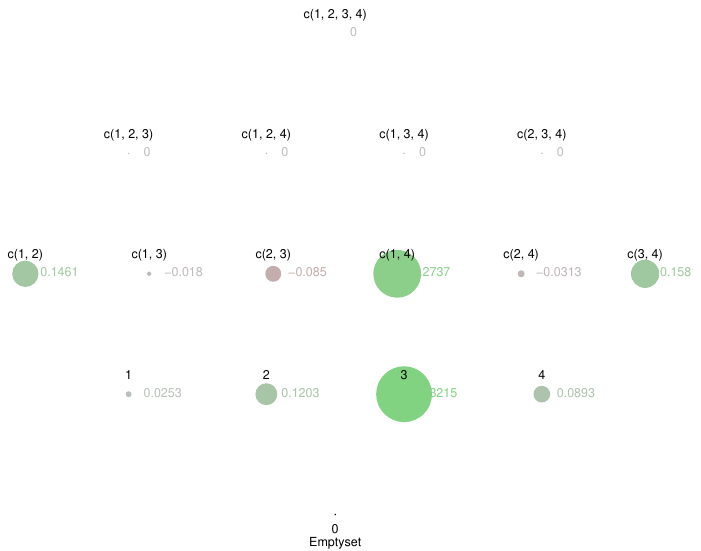}
		\end{minipage}%
	}%
	\subfigure[2-additive measure $\mu_4$  on $\{1,2,3,4,5\}$]{
		\begin{minipage}[t]{0.7\linewidth}
			\centering
			\includegraphics[width=3in]{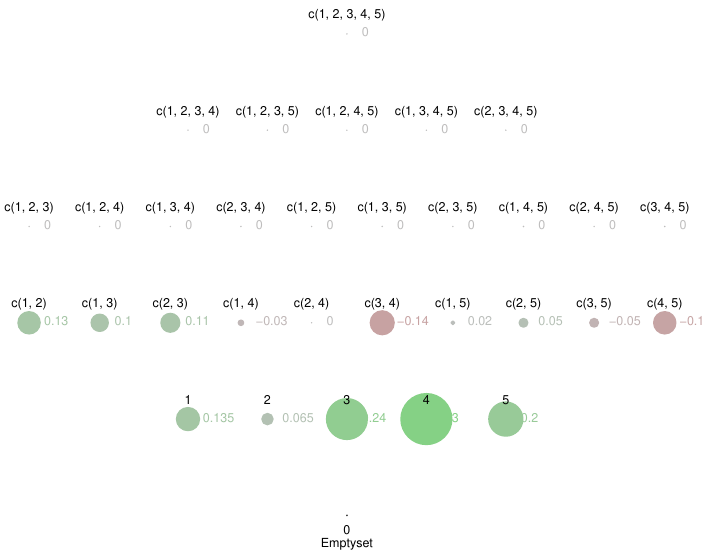}
		\end{minipage}%
	}%
	\caption{ k-additive measures with M\"obius representation indices.}
	\label{figs-2-additive-capa}
\end{figure}

A $k$-additive measure $\mu$ on $N$ only needs to consider the $\sum_{i=1}^k \binom{n}{i}$ Möbius representation values of subsets with cardinalities not exceeding $k$. A 1-additive measure is equivalent to an additive measure. As previously mentioned, a $k$-additive measure tends to accentuate the distinctions among subsets in the lower $k$ levels by setting the upper coefficients to zero, as illustrated in Figure \ref{figs-2-additive-capa}.

%
%

%



\begin{definition} \cite{marichal2007k}
Let $k \in \{1,..., n\}$, a fuzzy measure $\mu$ on $N$ is said to be $k$-tolerant if $\mu(A) = 1$ for all $A \subseteq N$ such that $|A| \geq k$ and there exists a subset
$ B \subseteq N$, with $|B| = k - 1$, such that $\mu(B) \neq 1$. A fuzzy measure $\mu$ on $N$ is said to be $k$-intolerant if $\mu(A) = 0$ for all $A \subseteq N$ such that $|A| \leq n - k$ and there exists a subset $ B \subseteq N$, with $|B| = n - k + 1$, such that $\mu(B) \neq 0$.
\end{definition}

\begin{definition} \cite{mesiar2018k,wuBeliakovkminitive}
	A fuzzy measure $\mu$ on $N$ is said to be $k$-maxitive if and only if $	\mu (A)=\bigvee_{B\subset A} \mu(B),$ $|A|\geq k+1$. A fuzzy measure $\mu$ on $N$ is said to be $k$-minitive if and only if $\mu (A)=\bigwedge_{B\supset A} \mu(B)$, $|A|\leq n-k-1$.
\end{definition}
 
\begin{figure}[!htp]
		\color{black}
	\centering	
	\subfigure[The topological graph of $\mu_5$ ]{
		\begin{minipage}[t]{0.4\linewidth}
			\centering
			\includegraphics[width=3in]{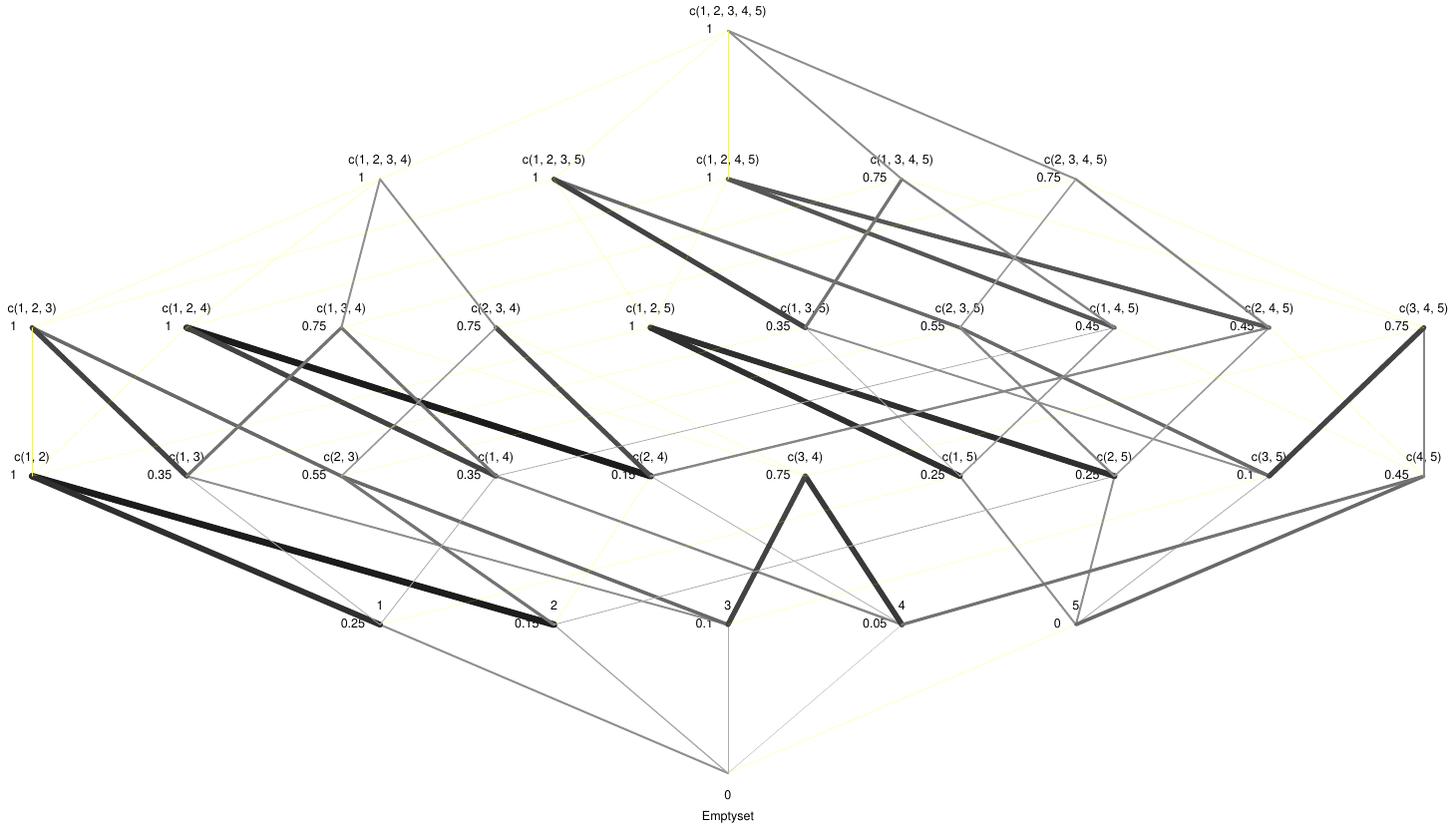}
		\end{minipage}%
	}%
	\subfigure[The topological  graph of $\mu_6$ ]{
		\begin{minipage}[t]{0.7\linewidth}
			\centering
			\includegraphics[width=3in]{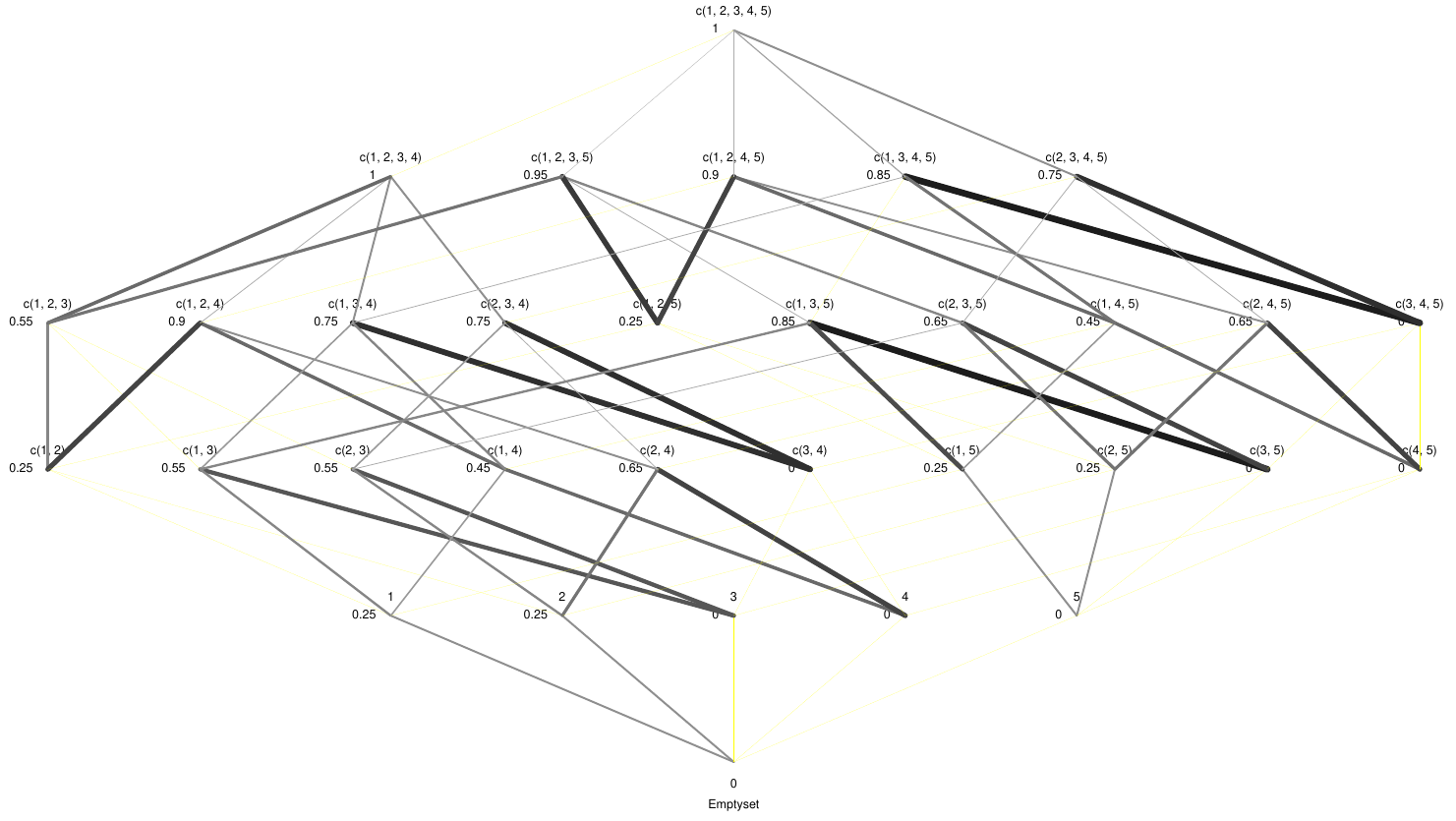}
		\end{minipage}%
	}%
		
		\subfigure[The "heigh-on" graph of $\mu_5$ ]{
			\begin{minipage}[t]{0.4\linewidth}
				\centering
				\includegraphics[width=3in]{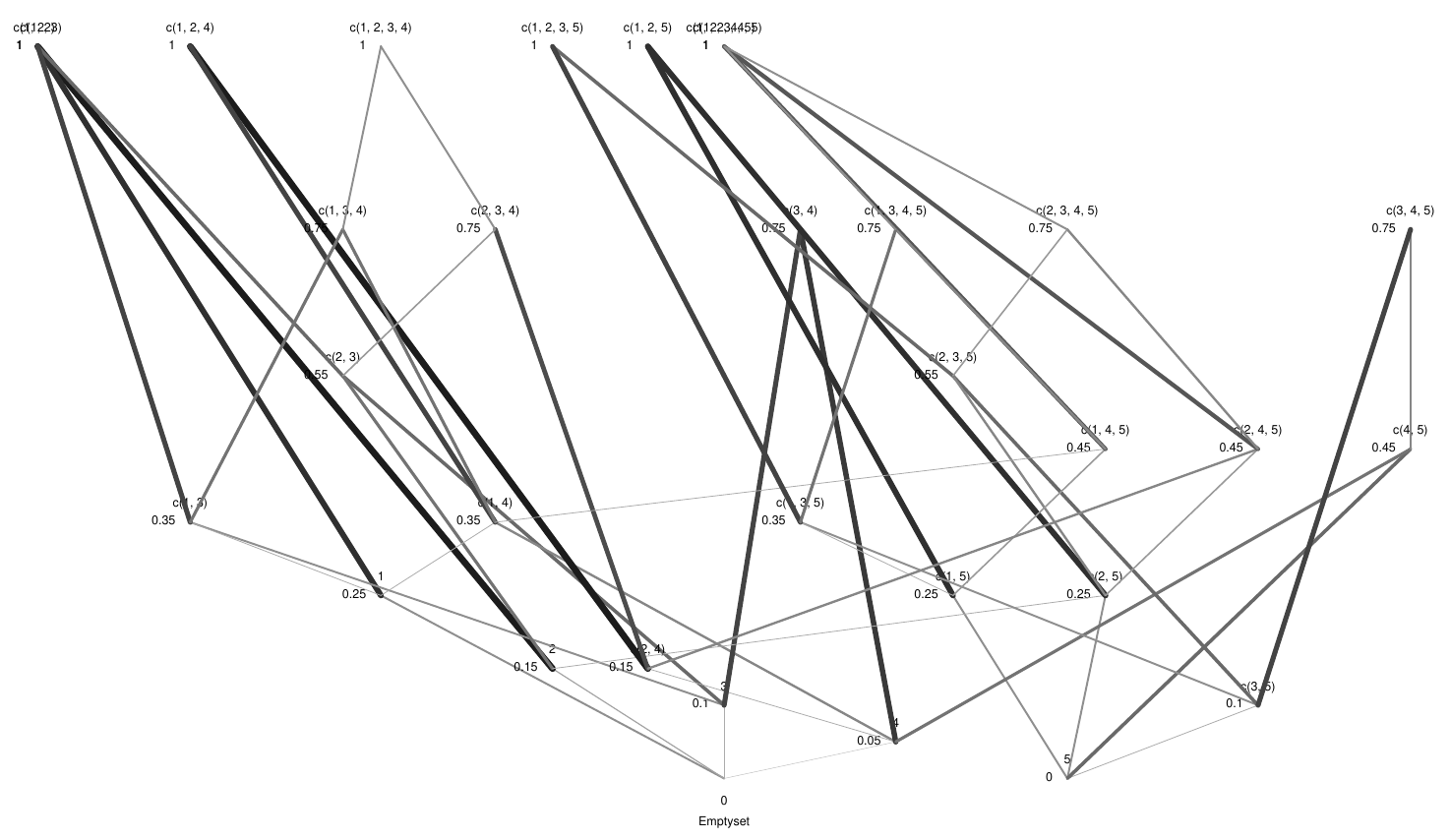}
			\end{minipage}%
		}%
		\subfigure[The "heigh-on" graph of $\mu_6$ ]{
			\begin{minipage}[t]{0.7\linewidth}
				\centering
				\includegraphics[width=3in]{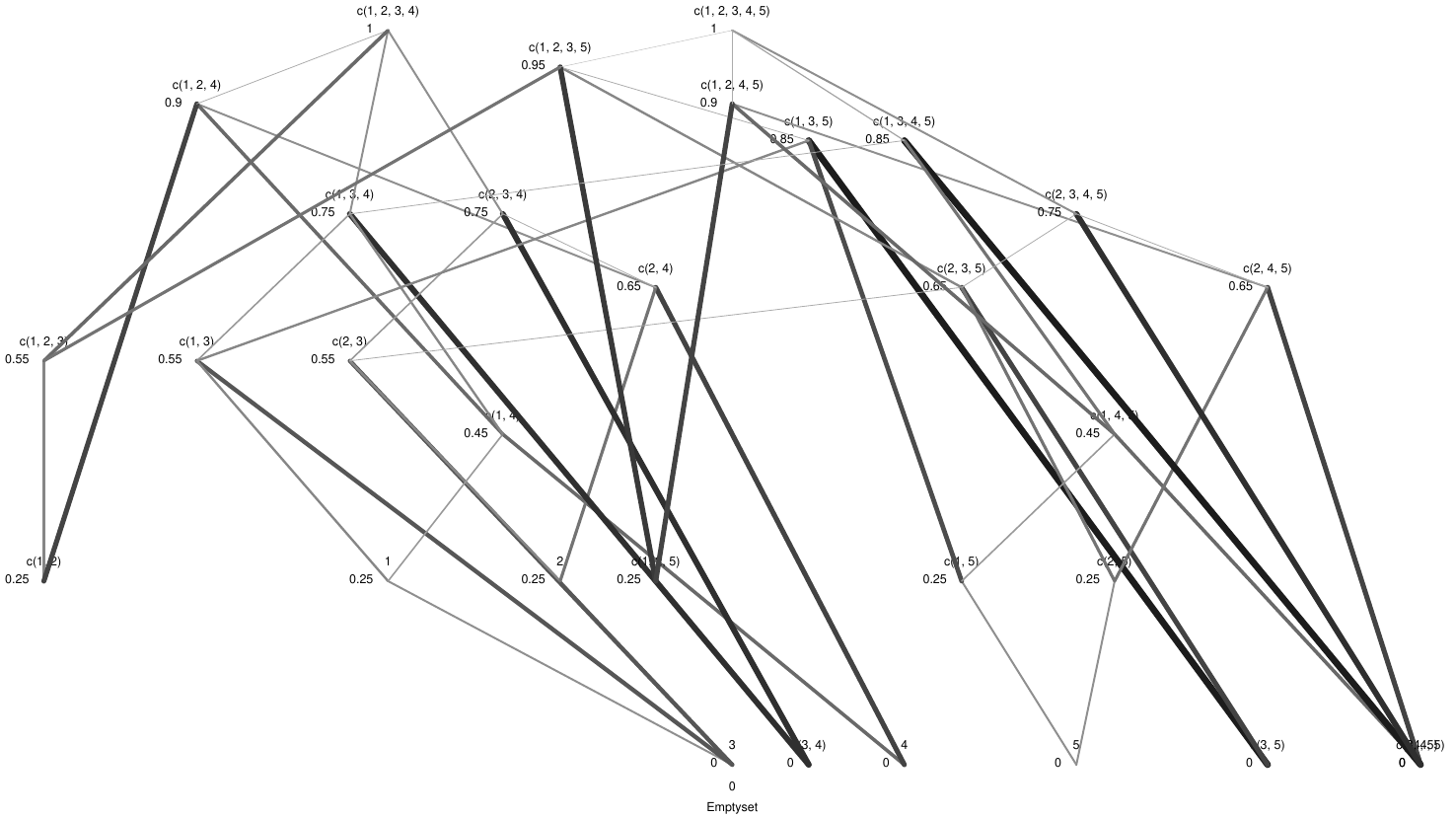}
			\end{minipage}%
		}%
		
		\centering
	\caption{ The 2-maxitive measure $\mu_3$ and 2-minitive measure $\mu_6$  on $\{1,2,3,4,5\}$}
	\label{fig-5-2-max-min-capa-new}
\end{figure}

It can be proved that $k$-tolerant and intolerant are special cases of the $k$-maxitive and minitive measures. One can see that for $k$-maximal fuzzy measure, the $k+1$ and higher levels' fuzzy measure values are just the maximum of those of its subsets and their heights will keep the largest heights of its subsets; to the contrary, the $n-k-s$ and lower levels of $k$-minitive measure will keep the smallest heights of its supersets, see many horizontal segments in sub figures (a) and (b) of Figure \ref{fig-5-2-max-min-capa-new}.

If a fuzzy measure exhibits multiple peaks, each with a height equal to 1, it is highly indicative of being a maximum fuzzy measure, with a single peak representing a normal fuzzy measure.
Conversely, if there are numerous bottom points, it suggests that the fuzzy measure is likely a minimum fuzzy measure, with a single bottom point corresponding to a normal fuzzy measure.
Additionally, if subsets at the $k+1$-level and higher levels reach peaks, it qualifies as a $k$-maxitive measure.
Conversely, if subsets at the $n-k+1$-level and lower levels reach bottom points, it falls into the category of a $k$-minitive measure.


An alternative approach to mitigate the complexity of fuzzy measures, as proposed in \cite{BeliakovWuLearning}, involves fixing the values of the fuzzy measure for all subsets with cardinalities exceeding $k$. This fixation is achieved by maximizing the partial entropy calculated over subsets with cardinalities greater than $k$.

\begin{definition}\label{def:kinter}
	A fuzzy measure  is called $k$-interactive if for some chosen $K\in [0,1]$ $$\mu(A)=K+\frac{a-k-1}{n-k-1}(1-K), \mbox{for all }A, a>k.$$
\end{definition}

\begin{figure}[!htp]
	\color{black}
	\centering	
	\subfigure[The "heigt-on" graph of $\mu_8$]{
		\begin{minipage}[t]{0.4\linewidth}
			\centering
			\includegraphics[width=3in]{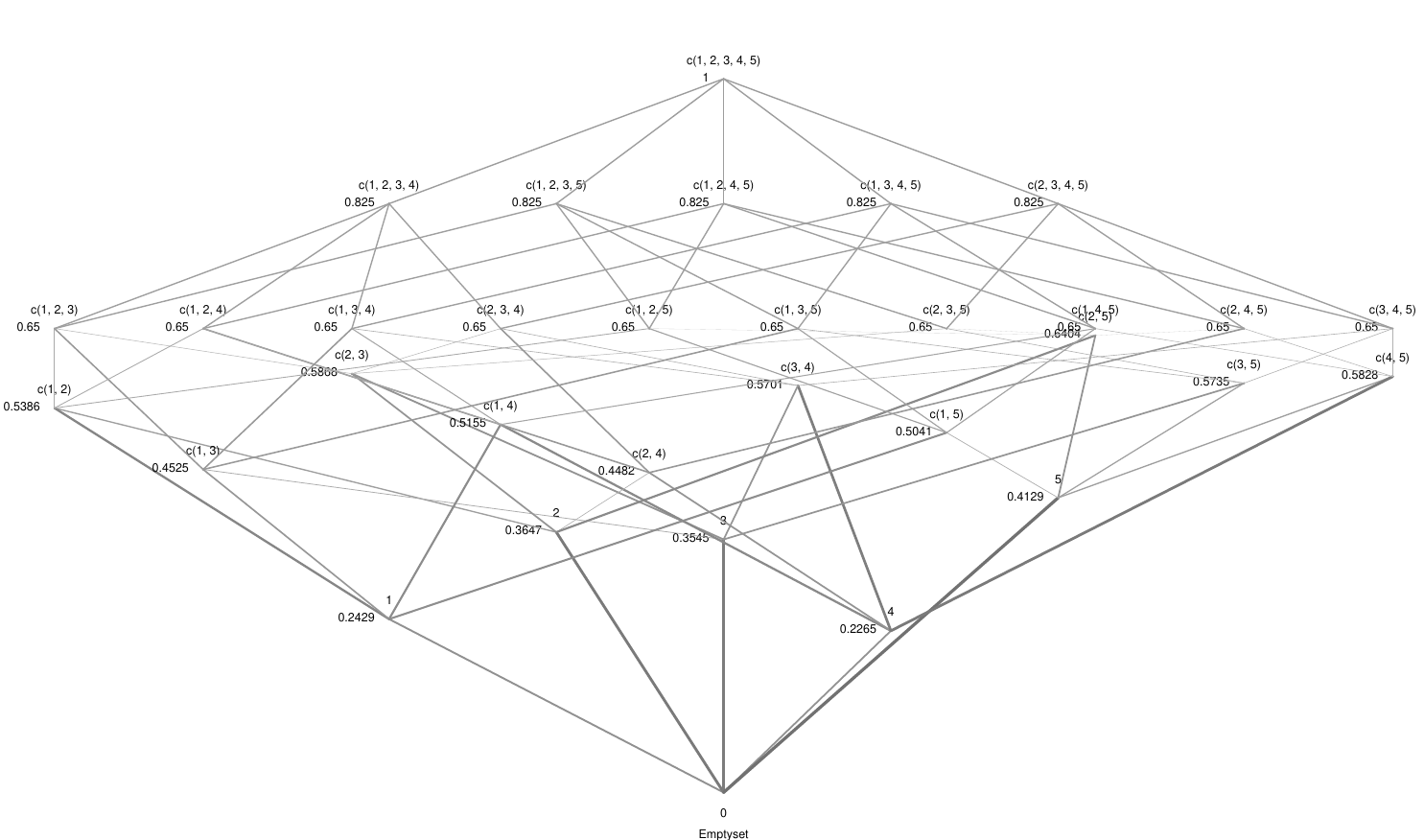}
		\end{minipage}%
	}%
	\subfigure[The "heigt-on" graph of $\mu_9$]{
		\begin{minipage}[t]{0.7\linewidth}
			\centering
			\includegraphics[width=3in]{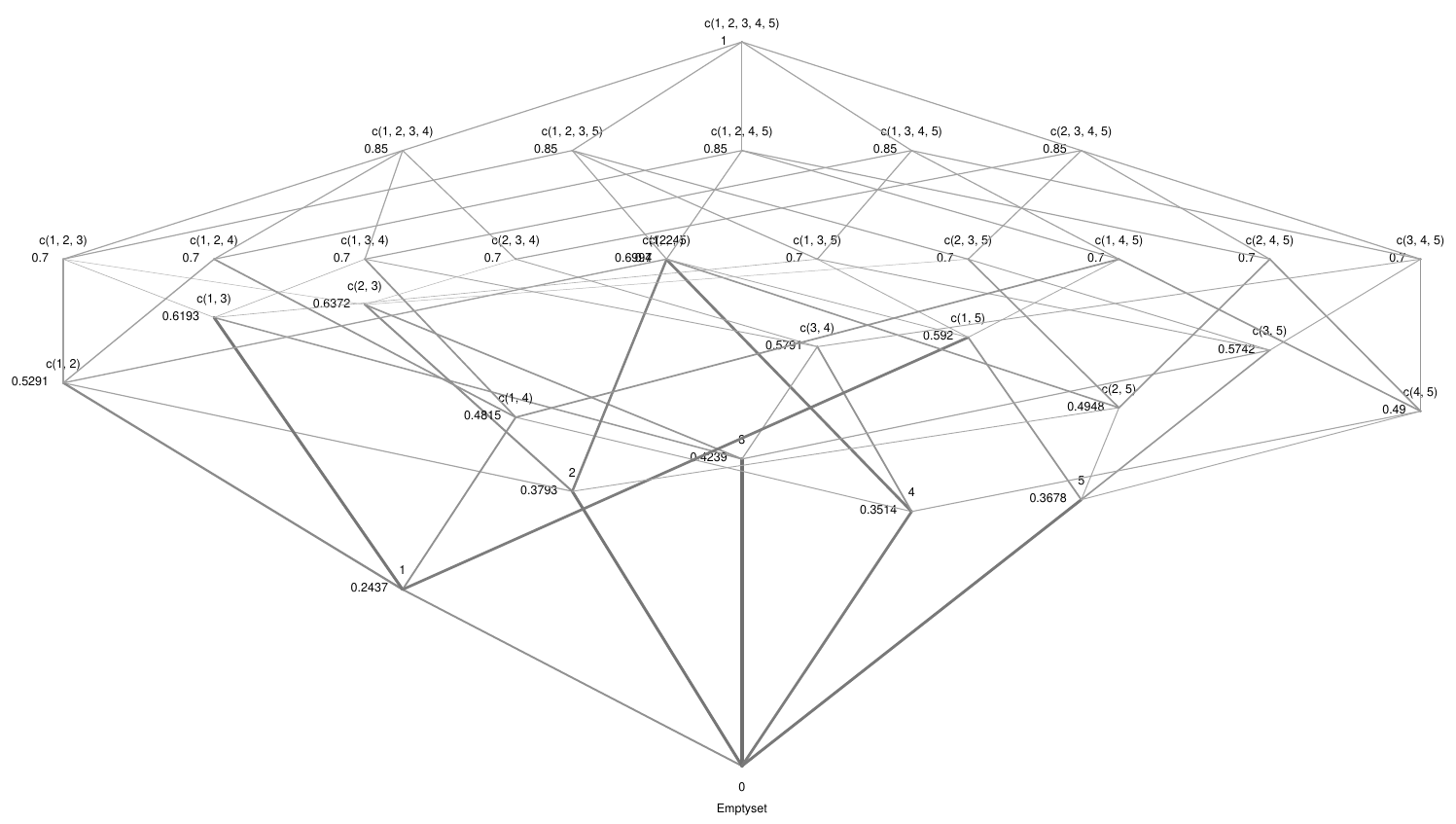}
		\end{minipage}%
	}%
	%
	\caption{Two 2-interactive measures $\mu_8$ and $\mu_9$ on $\{1,2,3,4,5\}$}
	\label{fig-5-criteria-2-interactive-capa-new}
\end{figure}

In this approach, the values $\mu(B) = K$ are fixed for all $B$ with $|B| = k+1$. In the representation of marginal contributions, all $\Delta_i(A) = \frac{1-K}{n-k-1}$ when $|A| > k$. These constraints serve to significantly reduce the number of variables and constraints associated with fuzzy measures. Notably, $k$-tolerant fuzzy measures emerge as special cases within this framework.
Interactions within subsets larger than $k$ are allowed but predefined by interactions within smaller subsets and the specified values of $k$ and $K$. The specific formula in Definition \ref{def:kinter} is designed to maximize the average contribution of the $n-k-1$ smallest inputs. Figure \ref{fig-5-criteria-2-interactive-capa-new} gives a 2-interactive measure.

\color{black}

\begin{figure}[!htp]
	\color{black}
	\centering
	\includegraphics[width=0.7\textwidth]{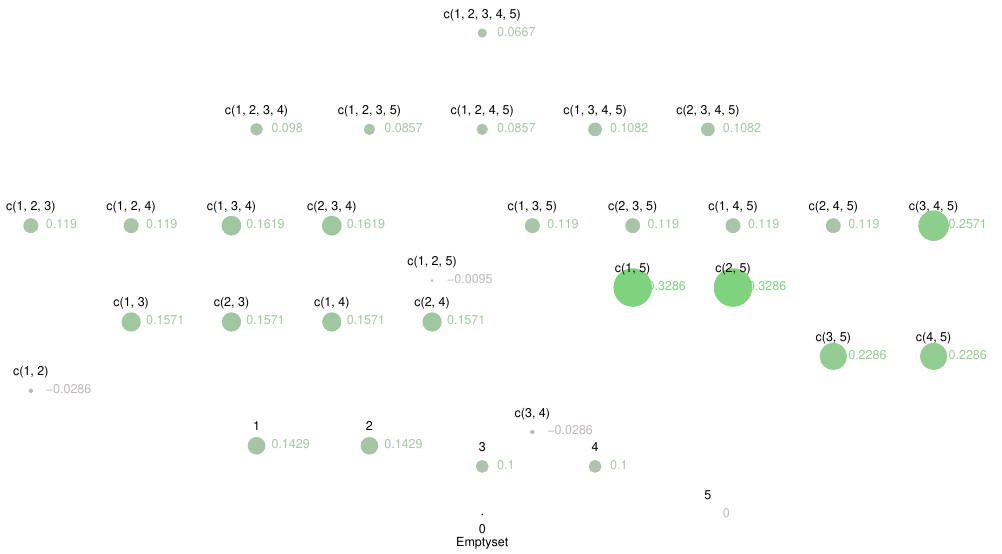}
	\caption{3-symmetric fuzzy measure $\mu_{10}$ on $\{1,2,3,4,5\}$ with basis $\{\{1,2\}, \{3,4\}, \{5\}\}$.} 
	\label{fig-5-3-symmetric-capa-new}
\end{figure}

\begin{definition} \cite{miranda2002p}
Let $\mu$ be a fuzzy measure on $N$, a subset $A \subseteq N$ is a subset of indifference with respect to $\mu$ if $\forall B_1,B_2 \subseteq A$, $|B_1| = |B_2|$, then $\mu(C\cup B_1) = \mu (C \cup B_2),\forall C \subseteq N\setminus A$.
%
A fuzzy measure $\mu $ on $N$ is said to be  $p$-symmetric if the coarsest partition of $N$ into subsets of indifference contains exactly $p$ subsets $A_1,..., A_p$, where $A_i$ is a subset of indifference, $A_i \cap A_j =\emptyset$,
$\cup_{i=1}^p A_i=N$, $i, j = 1,..., p$, and a partition $\pi$ is coarser than another partition $\pi'$
 if all subsets of $\pi$  are union of
some subsets of $\pi'$. The partition $\{A_1,...,A_p\}$ is called the basis of $\mu$.
\end{definition}
A 1-symmetric fuzzy measure is equivalent to a symmetric fuzzy measure.
When considering a basis $\{A_1,...,A_p\}$ for a $p$-symmetric fuzzy measure $\mu$ on $N$, each subset $S \subseteq N$ can be represented as a $p$-dimensional vector $\mathbf{b}_S=(b_1,..., b_p)$, where $b_i = |S \cap A_i|$, $i =1, ..., p$. In other words, a $p$-symmetric fuzzy measure necessitates the definition of $\prod_{i=1}^p(|A_i| + 1)$ coefficients \cite{miranda2002p}.
In Figure \ref{fig-5-3-symmetric-capa-new}, which illustrates the nonadditivity indices of all subsets with a basis of $\{\{1,2\}, \{3,4\}, \{5\}\}$, one can observe that the 3-symmetric fuzzy measure exhibits several groups of points with the same height.


\section{Nonlinear integral and its graphic representation}

Nonlinear integrals include many types of integrals \cite{beliakov2019discrete_book,klement2010universal,wang2010generalized}, the most common types are the Choquet integral\cite{choquet1954}, Sugeno integral \cite{sugeno1974theory} and Pan-integral\cite{wang2010generalized}.


\subsection{Choquet integral}

\begin{definition} \cite{choquet1954}
For a given $\mathbf{x}\in[0, + \infty]^n$, the discrete Choquet integral $\mathcal{C}$ of $\mathbf{x}$ with respect to fuzzy measure $\mu$ on $N$ is defined as follows:
\begin{equation*}\label{Choquet integral}
\mathcal{C}(\mathbf{x})=\sum_{i=1}^{n}(x_{(i)}-x_{(i-1)})\mu(\{(i),\ldots,(n)\}),
\end{equation*}
where $x_{(.)}$ is a non-decreasing permutation induced by $x_{i}$, $i=1, \ldots, n$,  i.e., $x_{(1)}\leq \ldots \leq x_{(n)}$, and  $x_{(0)}=0$ by convention.
\end{definition}

One way of calculating Choquet integral without ordering the evaluations in advance is given as: \cite{beliakov2019discrete_book}:
\begin{equation*}\label{eq:choquetbasis1}
\mathcal{C}_{\mu}(\mathbf x)=\sum_{ A \subseteq  N} \mu( A) \hat{x}_{ A}.
\end{equation*}
where basis function $
\hat{x}_{ A}=\max(0,\min_{i \in 
A}x_i-\max_{i \in  N \backslash A}x_i),\forall A \subseteq N.$

\begin{figure}[!htp]
	\centering
	\includegraphics[width=0.7\textwidth]{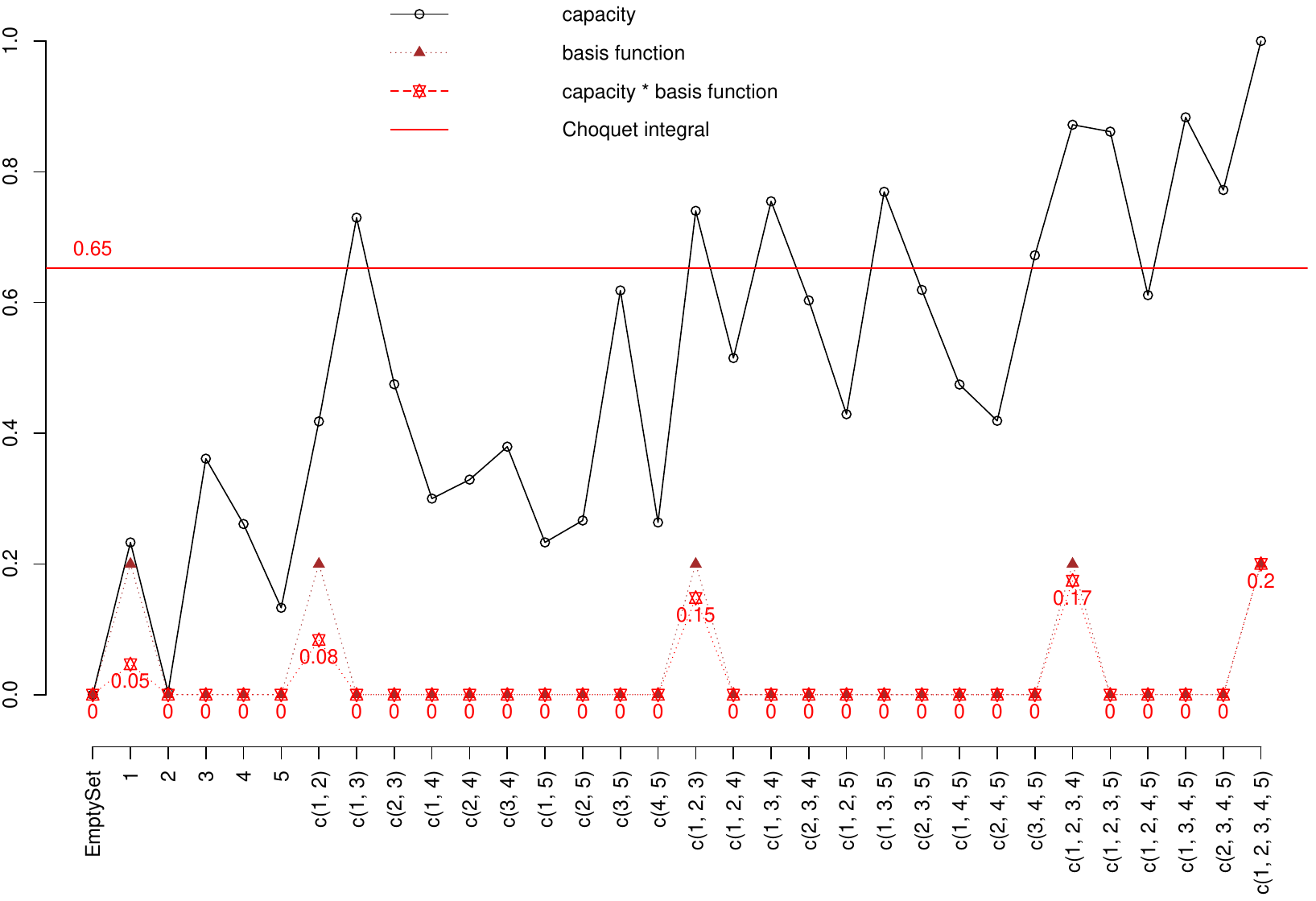}
	\caption{The Choquet integrals of (1,0.8,0.6,0.4,0.2).} 
	\label{fig-choquet-integral}
\end{figure}

\subsection{Sugeno integral}
\begin{definition} \cite{sugeno1974theory}
For a given $\mathbf{x}\in[0,1]^n$, the discrete Sugeno integral $\mathcal{S}$ of $x$ with respect to fuzzy measure $\mu$ on $N$ is defined as follows:
\begin{equation*}
\mathcal{S}(\mathbf{x})=\bigvee_{i=1}^{n}(x_{(i)}\wedge\mu(\{(i),\ldots,(n)\})),
\end{equation*}
where $x_{(.)}$ is a non-decreasing permutation induced by $x_{i}$, $i=1, \ldots, n$,  i.e., $x_{(1)}\leq \ldots \leq x_{(n)}$.
\end{definition}
From the above definition, we can get that \cite{dubois2001use} 
\begin{equation*}
\mathcal{S}_\mu(\mathbf{x})=\bigvee_{A \subseteq N} (\mu(A) \wedge \overline{x}_A),
\end{equation*}
where $\overline{x}_A= min_{i\in A}(x_i)$. 

\begin{figure}[!htp]
	\centering \color{black}
	\color{black}
	\includegraphics[width=0.6\textwidth]{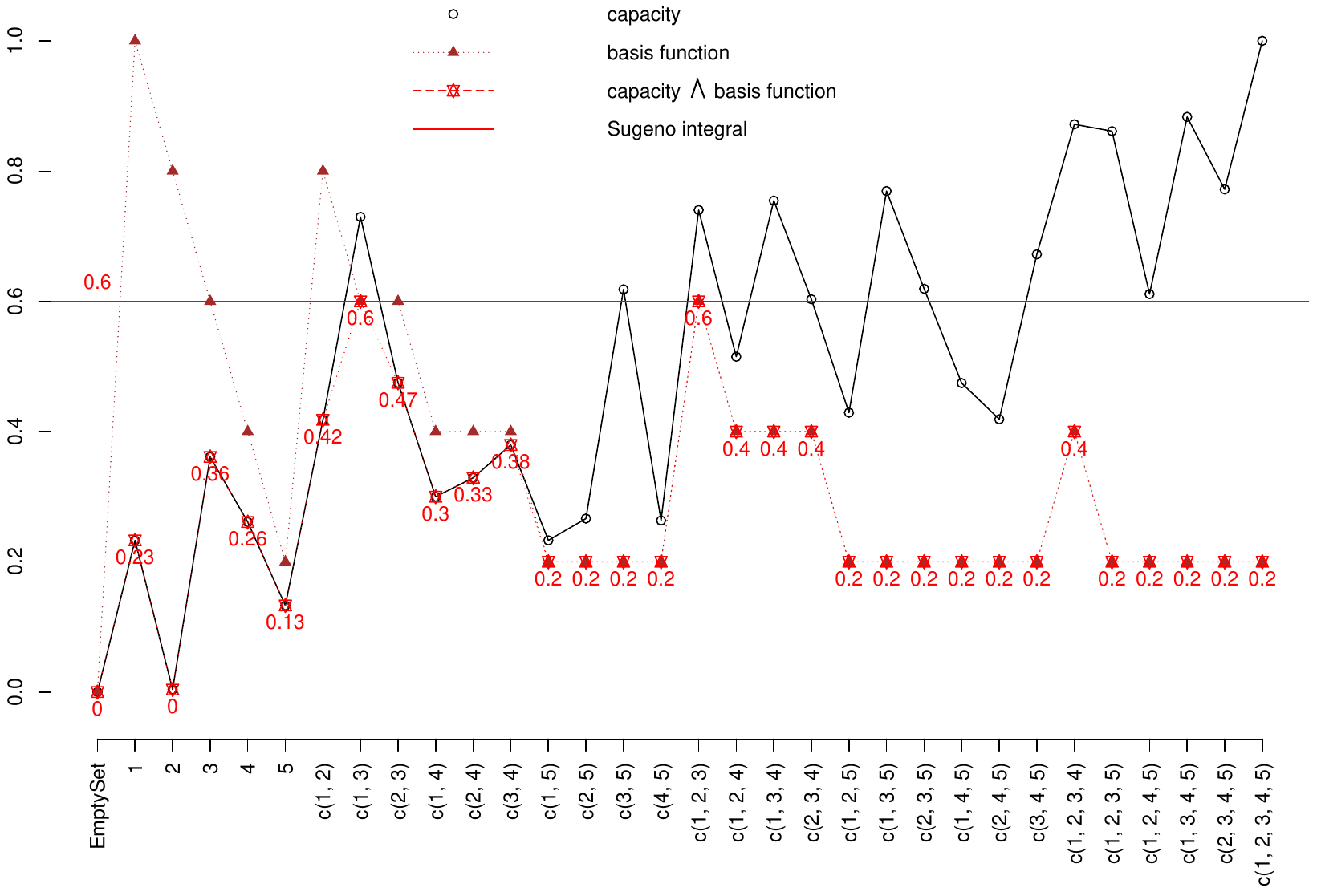}
	\caption{The Sugeno integrals of (1,0.8,0.6,0.4,0.2).} 
	\label{fig-sugeno-integral}
\end{figure}


\subsection{Pan-integral}
\begin{definition} \cite{wang2010generalized}
For a given $\mathbf{x}\in[0,1]^n$, the discrete pan integral $\mathcal{N}$ of $x$ with respect to fuzzy measure $\mu$ on $N$ is defined as follows:
\begin{equation*}
\mathcal{N}_\mu (\mathbf{x})=\bigvee_{i=1}^{n}(x_{(i)}\mu(\{(i),\ldots,(n)\})),
\end{equation*}
where $x_{(.)}$ is a non-decreasing permutation induced by $x_{i}$, $i=1, \ldots, n$,  i.e., $x_{(1)}\leq \ldots \leq x_{(n)}$.
\end{definition}
 With the same basis function of Sugeno integral of fuzzy measure form,  the pan integral can be represented as \cite{beliakov2019discrete_book}:
\begin{equation*}
\mathcal{N}_\mu(\mathbf{x})=\bigvee_{A \subseteq N} (\mu(A) \overline{x}_A),
\end{equation*}
where $\overline{x}_A= min_{i\in A}(x_i)$. From the non preordered representations of Sugeno and pan integrals, we have, for same $\mu$ and  $\mathbf{x}$, $\mathcal{N}_\mu(\mathbf{x}) \le \mathcal{S}_\mu(\mathbf{x})$.

\begin{figure}[!htp]
	\centering \color{black}
	\includegraphics[width=0.6\textwidth]{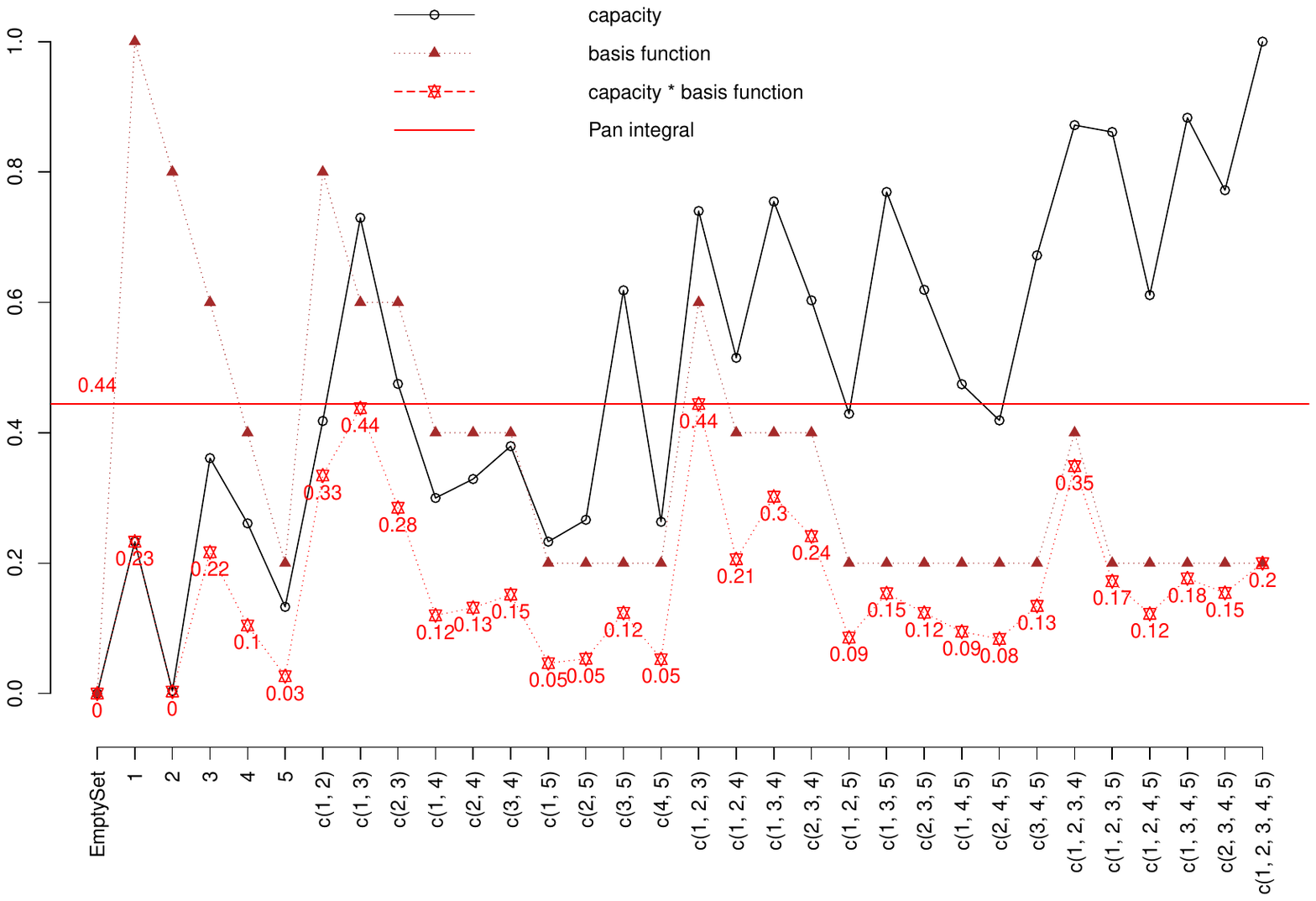}
	\caption{The Pan integrals of (1,0.8,0.6,0.4,0.2).} 
	\label{fig-pan-integral}
\end{figure}

\begin{figure}[!htp]
	\centering
	\includegraphics[width=0.9\textwidth]{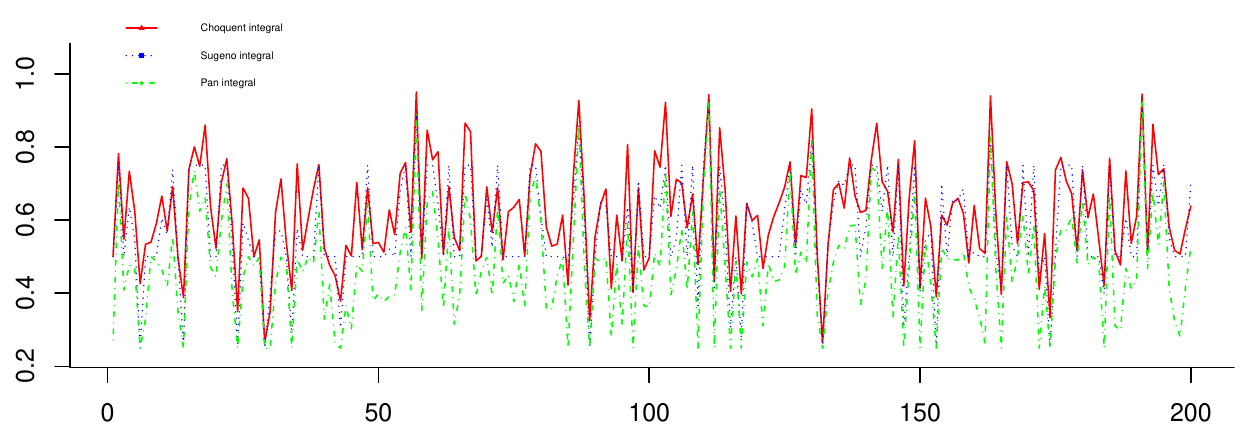}
	\caption{The integrals of (0.2,0.5,0.75,1) w.r.t. random generated fuzzy measures.} 
	\label{fig-3-integral}
\end{figure}

The figure \ref{fig-3-integral} shows an empirical comparison of these three integrals of alternative (0.2,0.5,0.75,1) with 200 random generated fuzzy measures, see \cite{combarro2006identification} for the detailed algorithm. One can conclude that with high probability three integrals' values are likely to hold the relationship: Choquet integral > Sugeno integral > pan integral.


\section{Graphic representation of fuzzy measure fitting}

The fuzzy measure fitting or identifying is one of the core focuses for fuzzy measure theory and applications. Many methods have been proposed and studied, such as least-squares principle method \cite{grabisch2008review,grabisch2010decade,Meyer2005Choice}, the least absolute deviation criterion method \cite{beliakov2009construction}, the maximum split method \cite{marichal2000determination,Meyer2005Choice}, the maximum entropy principle method \cite{Kojadinovic2007minimum,kojadinovic2005axiomatic,marichal2002entropy}, the maximum log-likelihood principle method
\cite{Fallah2012Learning,Hullermeier2013Efficient,Fallah2010preference},  the compromise principle method \cite{wu2014compromise} and MCCPI (multicriteria correlation preference information) based minimum divergence principle method \cite{wu20152}.

However, these fuzzy measure fitting methods generally only provide the final optimal solutions and seldom reveal the fitting process. By virtue of graphic representation of fuzzy measure, it becomes relatively easy to show the identification process gradually. 

Taking the linear programming based least absolute deviation criterion method as an example, its fitting model is given as follows:
\begin{equation*} \label{eq-lp-absolute-deviation-model}
\begin{split}
&\min z=\sum_{x \in L} d_x^+ +d_x^- \\
s.t. &\begin{cases}
\text{The boundary and monotonicity condition of fuzzy measure},\\
\text{The decision preference information constraints},\\
C(\textbf{x})-d_x^++d_x^-=y(x), x \in L,\\
\end{cases}
\end{split}
\end{equation*}
where $L$ is the set of history alternatives, $y(x)$ is the desired overall evaluation of $x$, $d_x^+, d_x^-$ are the positive and negative deviation variables.
We adapt the example in \cite{grabisch2008review} and suppose there are seven alternatives whose partial evaluations on five criteria and desired overall evaluations are shown in Table \ref{table-7-alter-scores}.

\begin{table}[!htb]
	\centering
		\caption{Scores of 7 alternatives on 5 criteria} 
	\label{table-7-alter-scores}
	\resizebox{\textwidth}{17mm}{
	\begin{tabular}{rrrrrrr}
		\hline
	Alternatives	& Criterion 1 & Criterion 2 & Criterion 3 & Criterion 4 & Criterion 5 & Overall evaluations \\ 
		\hline
		1 & 18.0 & 11.0 & 11.0 & 11.0 & 18.0 & 15.0 \\ 
		2 & 18.0 & 11.0 & 18.0 & 11.0 & 11.0 & 14.5 \\ 
		3 & 11.0 & 11.0 & 18.0 & 11.0 & 18.0 & 14.0 \\ 
		4 & 18.0 & 18.0 & 11.0 & 11.0 & 11.0 & 13.5 \\ 
		5 & 11.0 & 11.0 & 18.0 & 18.0 & 11.0 & 13.0 \\ 
		6 & 11.0 & 11.0 & 18.0 & 11.0 & 11.0 & 12.5 \\ 
		7 & 11.0 & 11.0 & 11.0 & 11.0 & 18.0 & 12.0 \\ 
		\hline
	\end{tabular}
}
\end{table}

\color{black}

\begin{figure}[!htp]
	\color{black}
	\centering	
	\subfigure[After having trained by alternative 1.]{
		\begin{minipage}[t]{0.5\linewidth}
			\centering
			\includegraphics[width=2.5in]{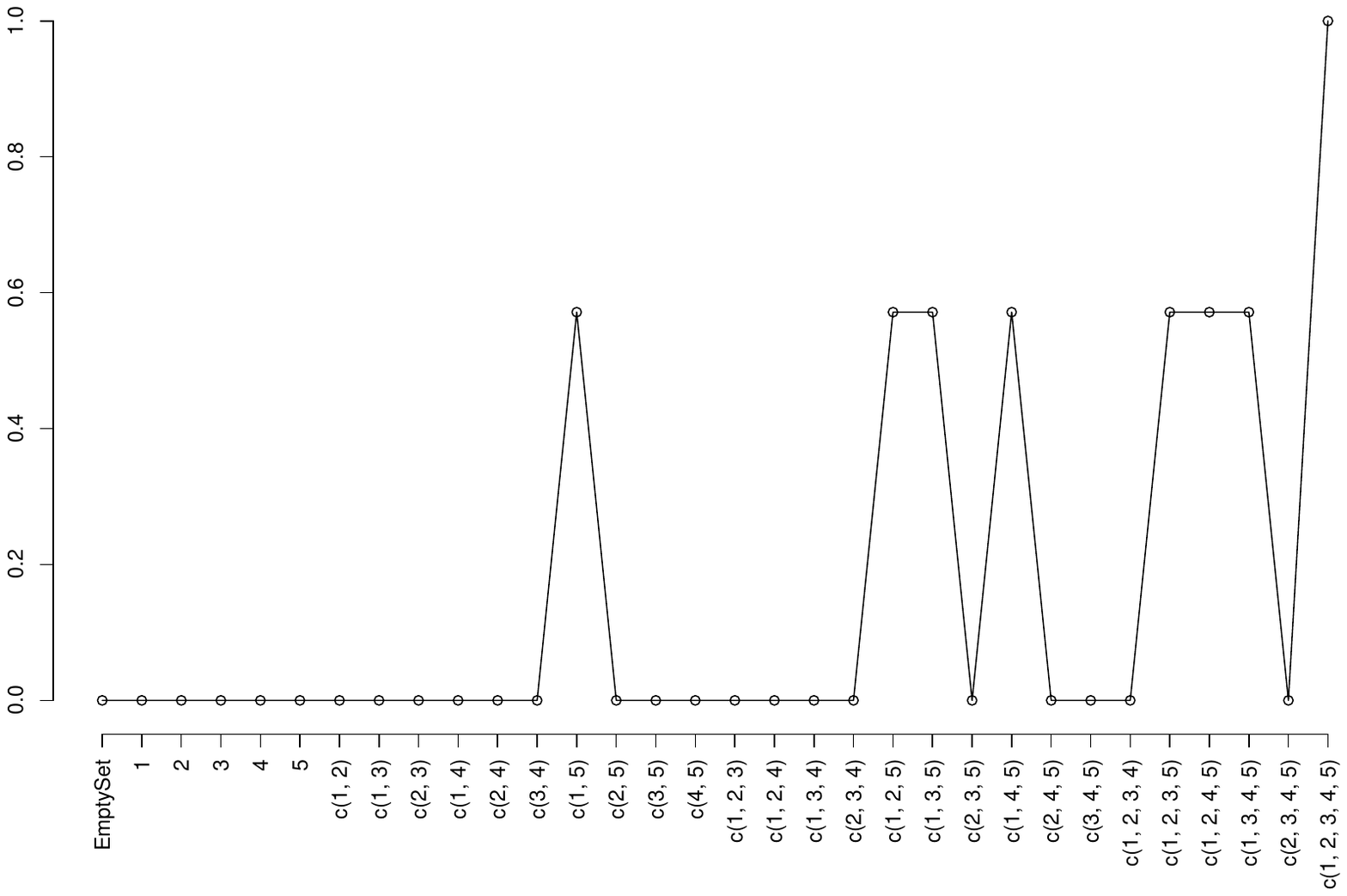}
		\end{minipage}%
	}%
	\subfigure[After having trained by alternatives 1 $\sim$ 2.]{
		\begin{minipage}[t]{0.5\linewidth}
			\centering
			\includegraphics[width=2.5in]{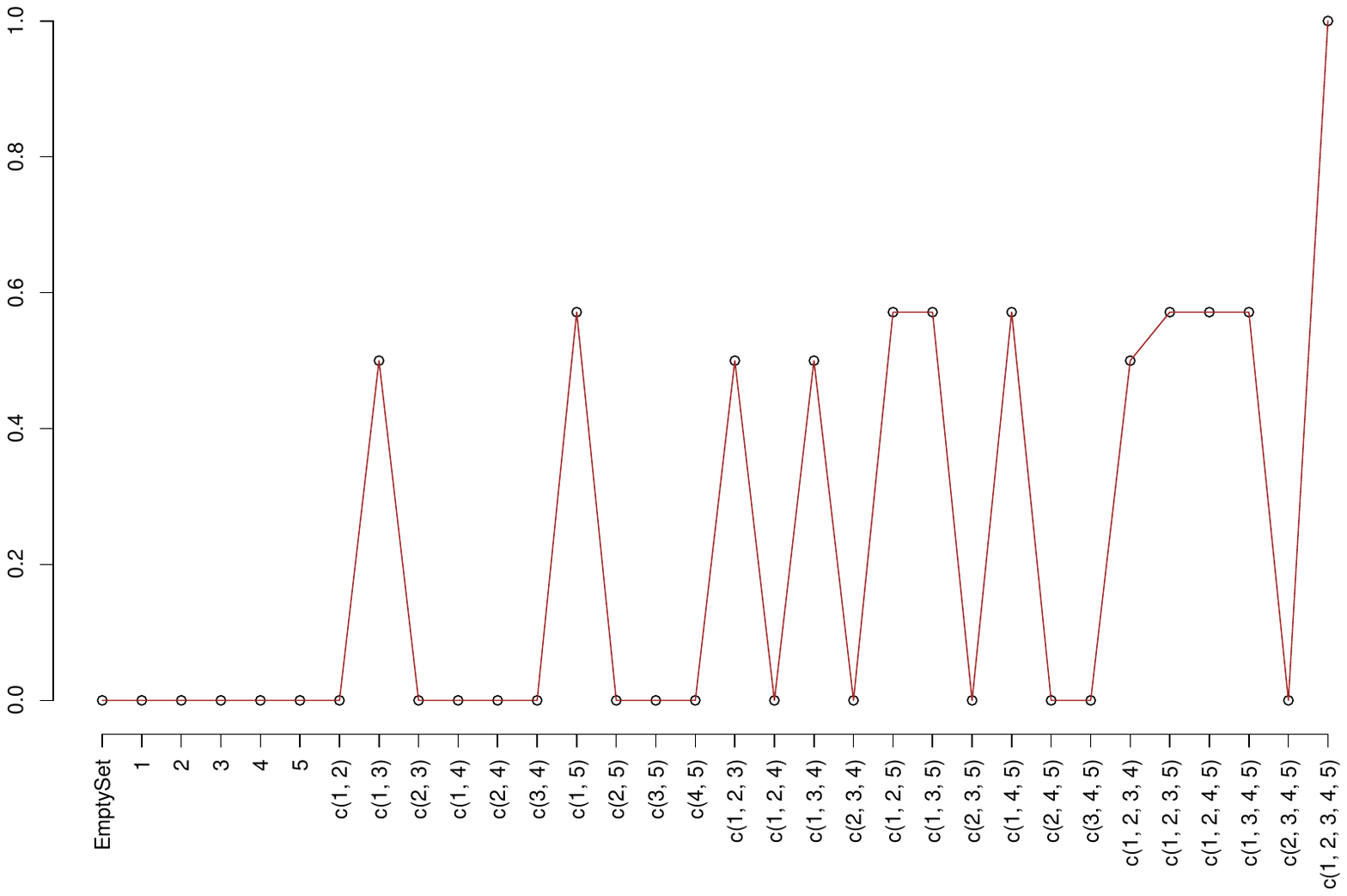}
		\end{minipage}%
	}
	
	\subfigure[After having trained by alternatives 1 $\sim$ 3.]{
		\begin{minipage}[t]{0.5\linewidth}
			\centering
			\includegraphics[width=2.5in]{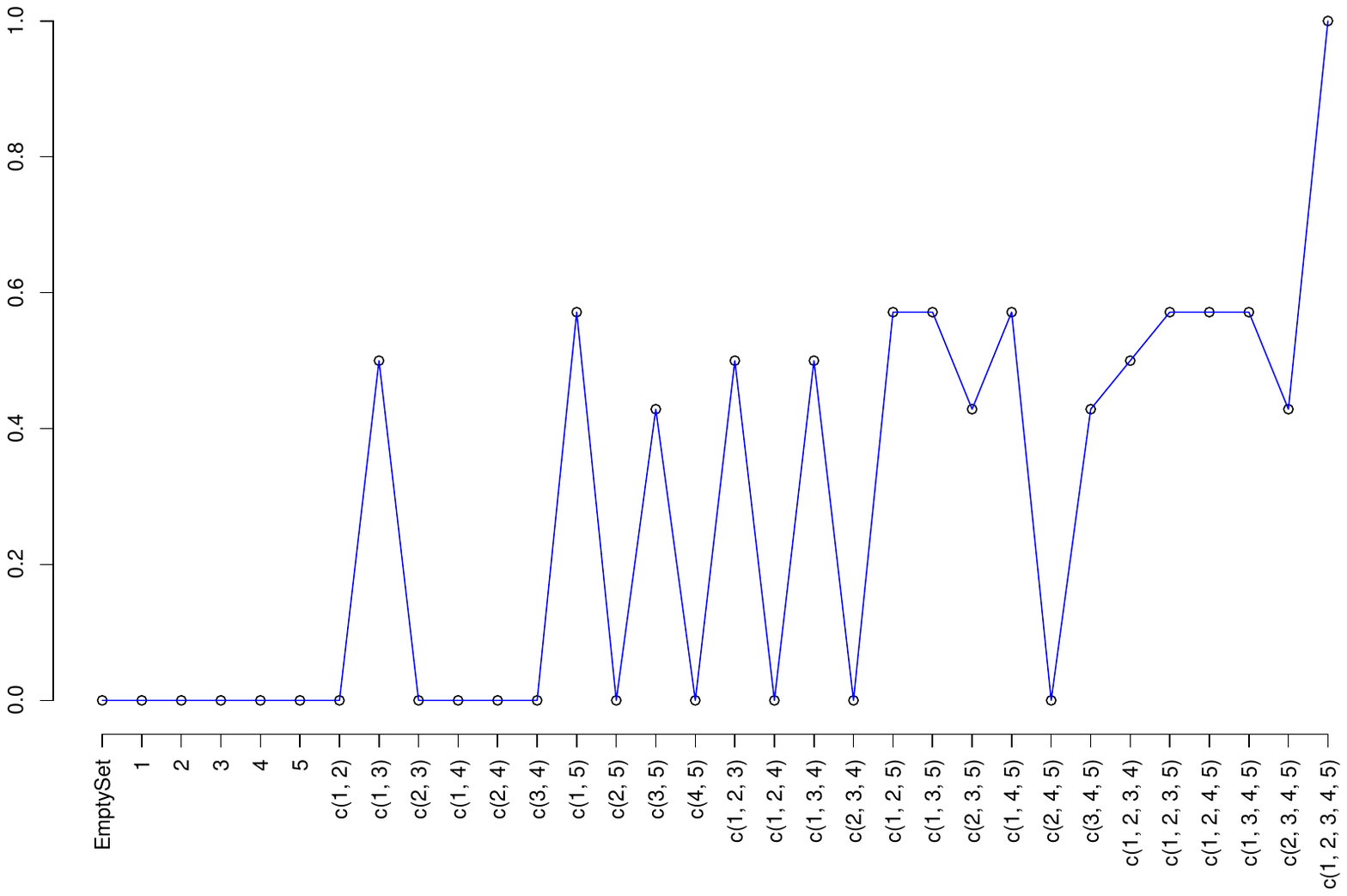}
		\end{minipage}%
	}%
	\subfigure[After having trained by alternatives 1 $\sim$ 4.]{
		\begin{minipage}[t]{0.5\linewidth}
			\centering
			\includegraphics[width=2.5in]{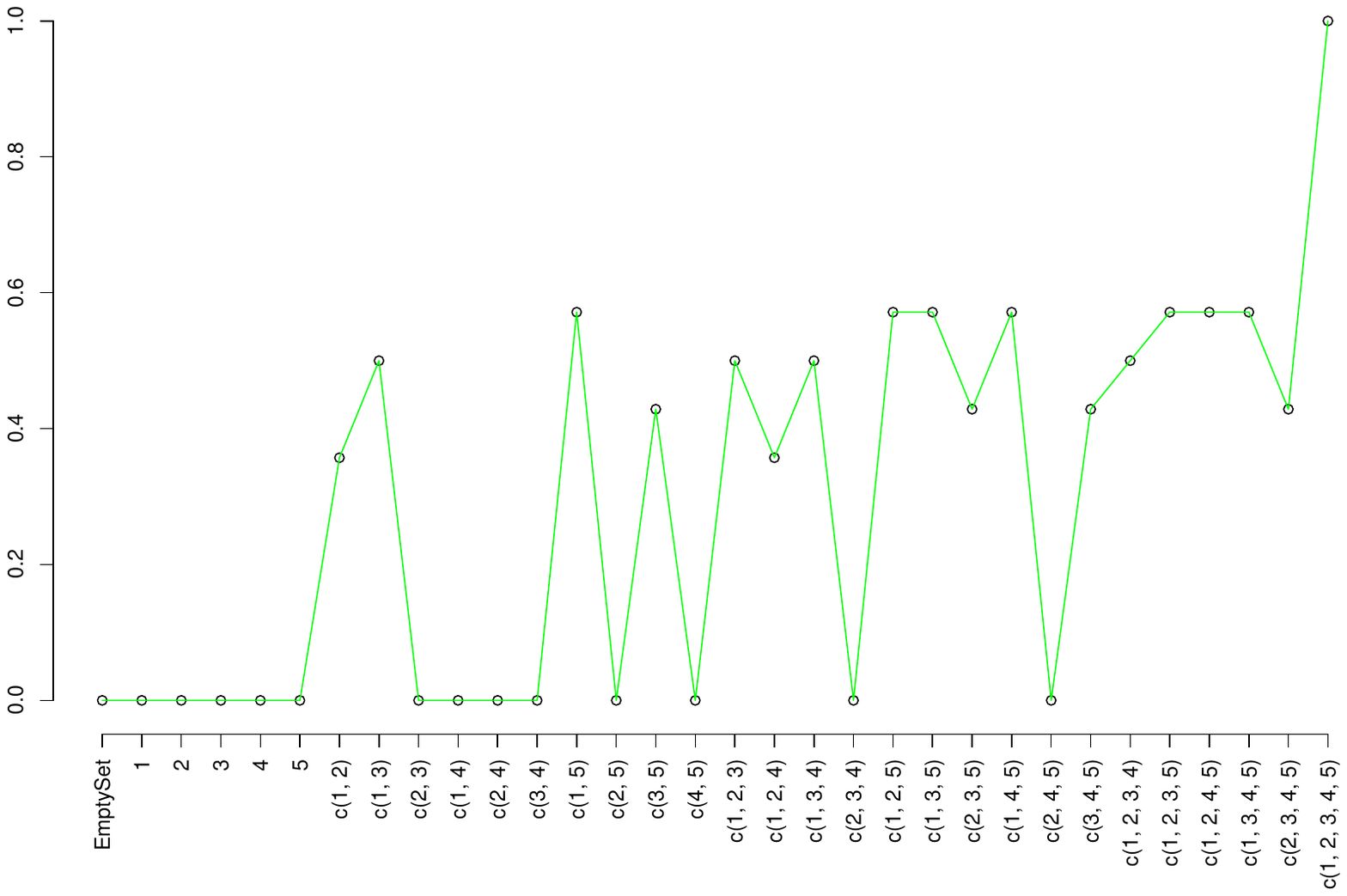}
		\end{minipage}%
	}
	\quad
	\subfigure[After having trained by alternatives 1 $\sim$ 5.]{
		\begin{minipage}[t]{0.5\linewidth}
			\centering
			\includegraphics[width=2.5in]{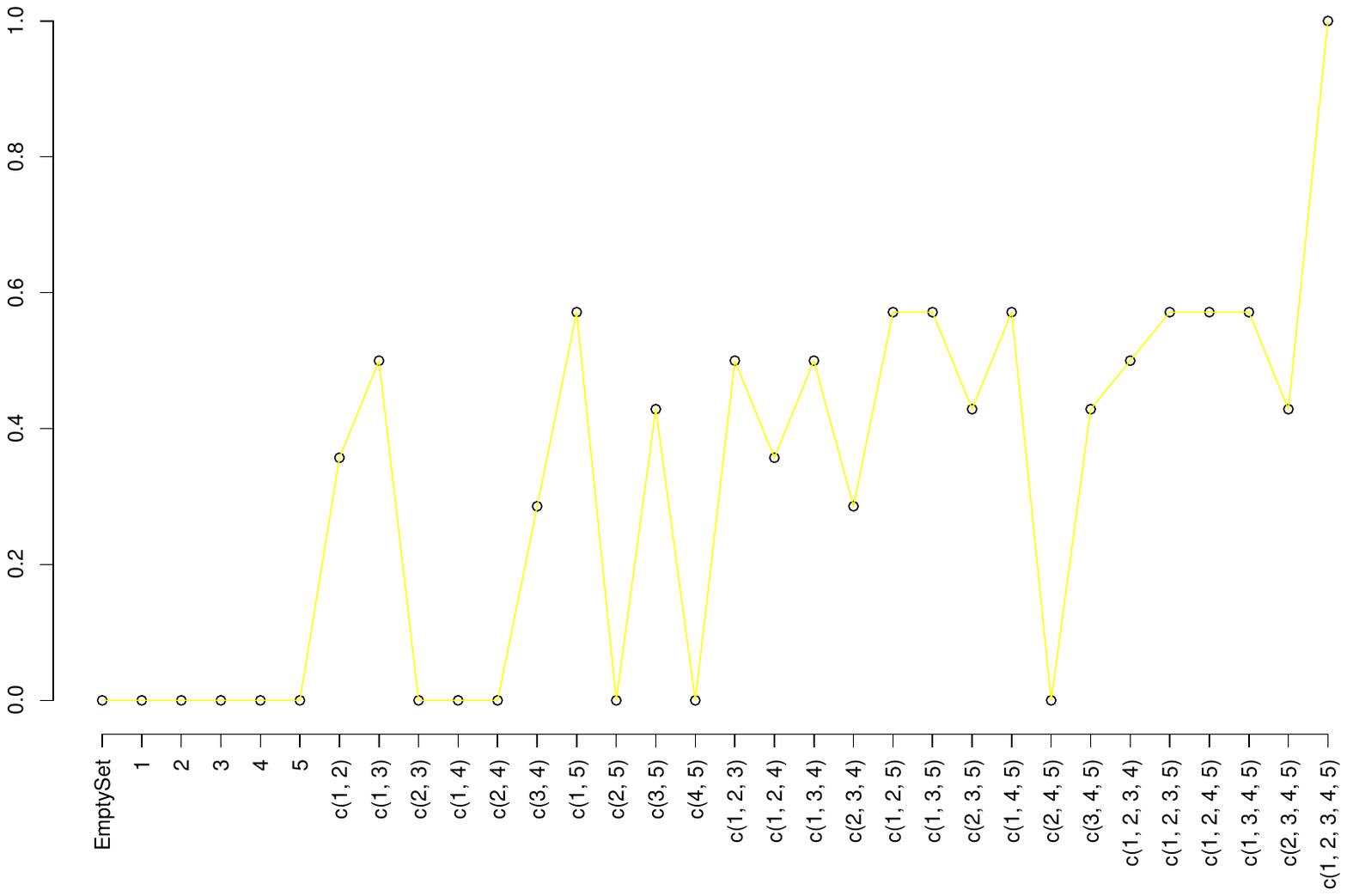}
		\end{minipage}%
	}%
	\subfigure[After having trained by alternatives 1 $\sim$ 6.]{
		\begin{minipage}[t]{0.5\linewidth}
			\centering
			\includegraphics[width=2.5in]{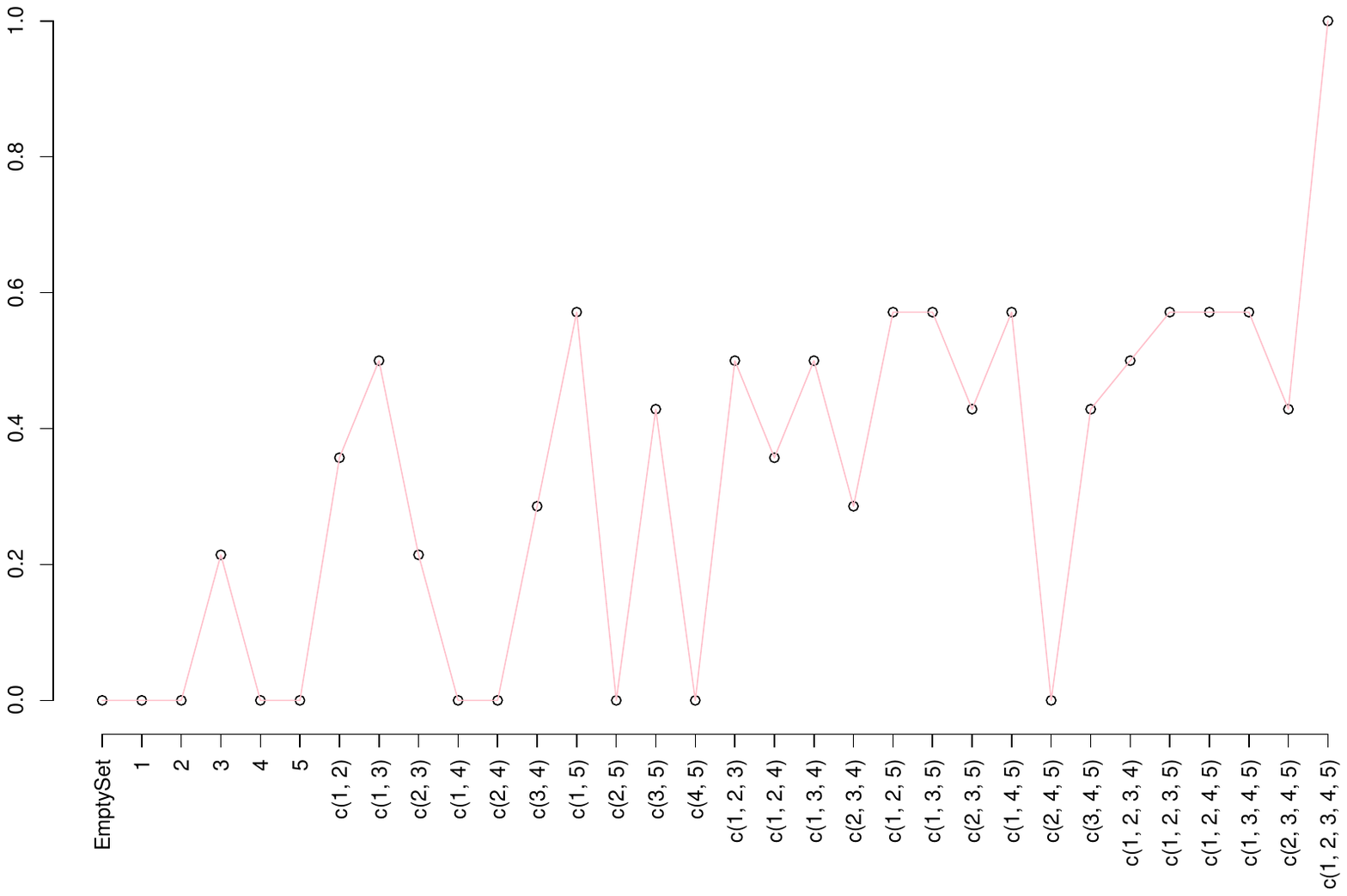}
		\end{minipage}%
	}
	
	\subfigure[After having trained by alternatives 1 $\sim$ 7.]{
		\begin{minipage}[t]{0.5\linewidth}
			\centering
			\includegraphics[width=2.5in]{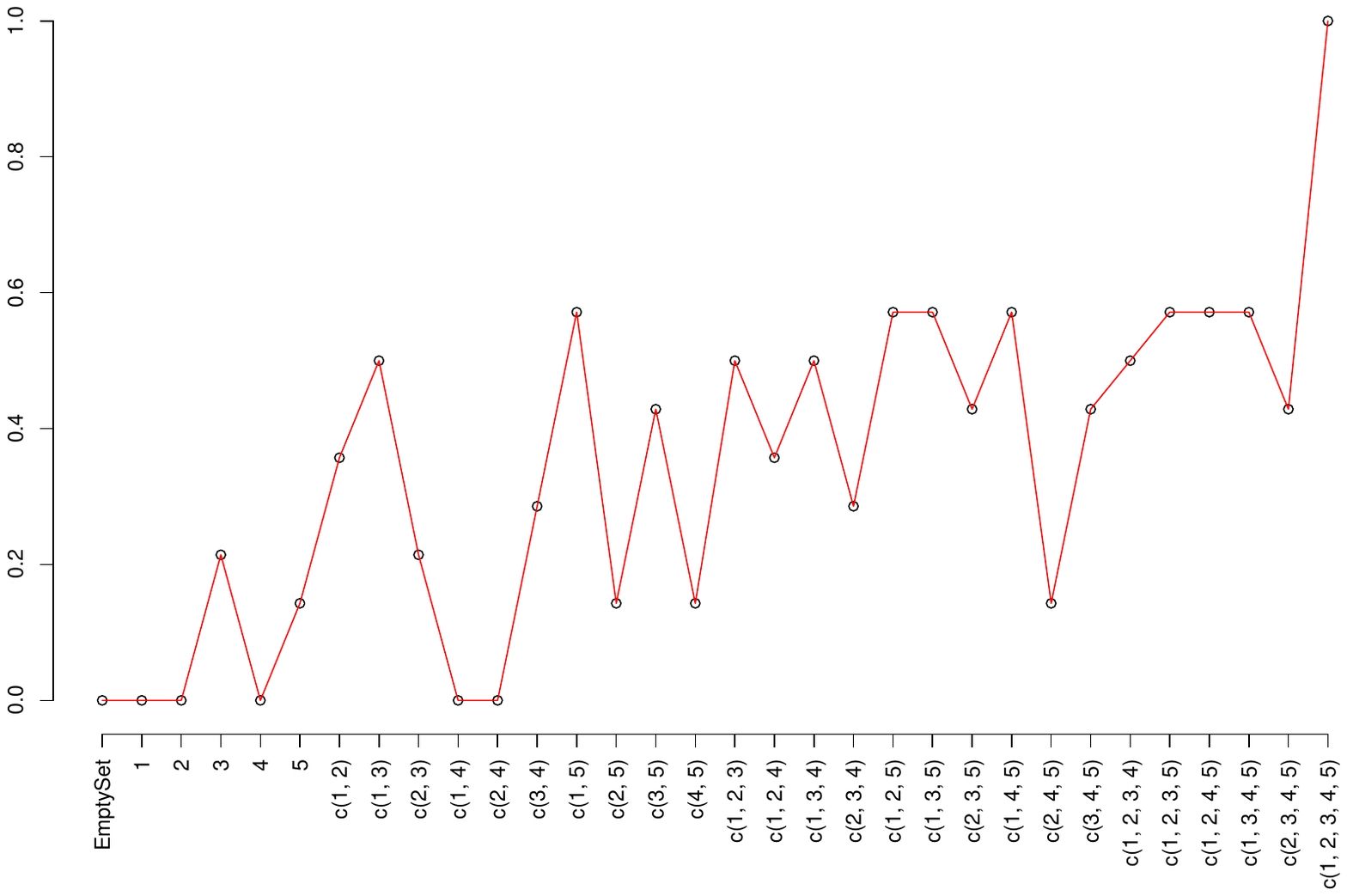}
		\end{minipage}%
	}%
	\caption{The process of fitting fuzzy measure $\mu_{11}$ by 7 alternatives.}
	\label{fig-7-capa-learning}
\end{figure}

\begin{figure}[!htp]
	\centering 	\color{black}
	\includegraphics[width=0.6\textwidth]{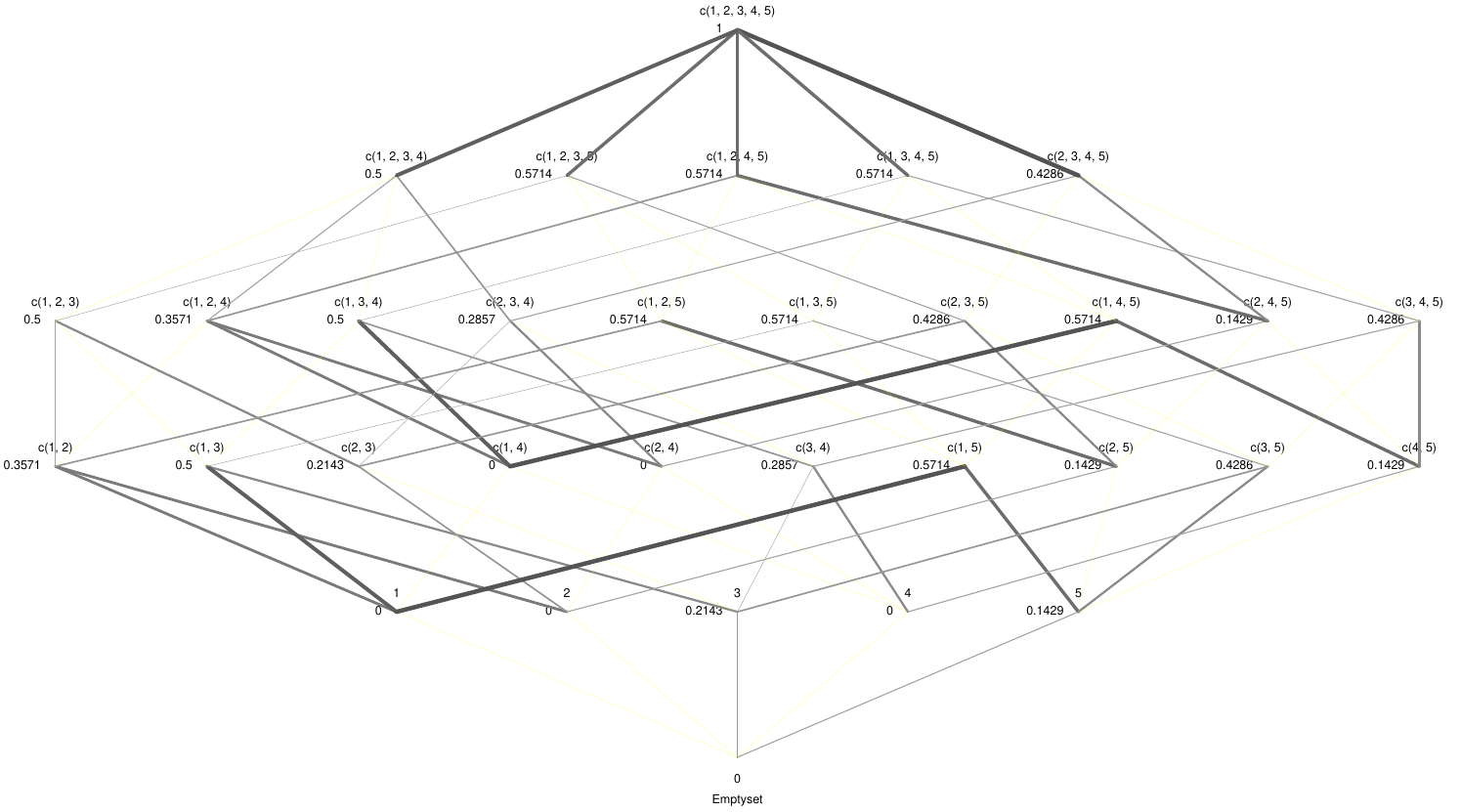}
	\caption{ The topology graph of fuzzy measure $\mu_{11}$} 
	\label{fig-fitted-cap-heigh-off}
\end{figure}

By executing the linear programming model with the seven alternatives in sequence, we can illustrate the changing or evolving process of these optimal fuzzy measures in graphical form, as shown in Figure \ref{fig-7-capa-learning}.
Initially, the minimum fuzzy measure is set, where the fuzzy measure value for the universal set is equal to 1, i.e., $\mu\{1,2,3,4,5\}=1$, and the fuzzy measure values for all other subsets are set to 0.
As we train the model with the first record in Table \ref{table-7-alter-scores}, the obtained fuzzy measure is displayed in subfigure (a), where $\mu\{1,5\}=0.5714$ (and correspondingly, its supersets, except the universal set, also have fuzzy measure values of 0.5714).
In the second round, shown in subfigure (b), we have $\mu\{1,3\}=0.5$. In the third round, depicted in subfigure (c), we further obtain $\mu\{3,5\}=0.4285$.
This process continues with each round, and the obtained fuzzy measures evolve as follows: in the fourth round (subfigure (d)), $\mu\{1,2\}=0.3571$; in the fifth round (subfigure (e)), $\mu\{3,4\}=0.2857$; in the sixth round (subfigure (f)), $\mu\{3\}=0.2142$; and finally, in the seventh and final round, we achieve $\mu\{5\}=0.1429$.
Figure \ref{fig-fitted-cap-heigh-off} provides a more detailed view of the final fuzzy measure's topological graph with "height-on," where the presence of many zero marginal contributions (represented by the pale yellow segments) indicates that certain aspects are still under training.


\begin{figure}[!htp]
	\color{black}
	\centering
	\includegraphics[width=0.8\textwidth]{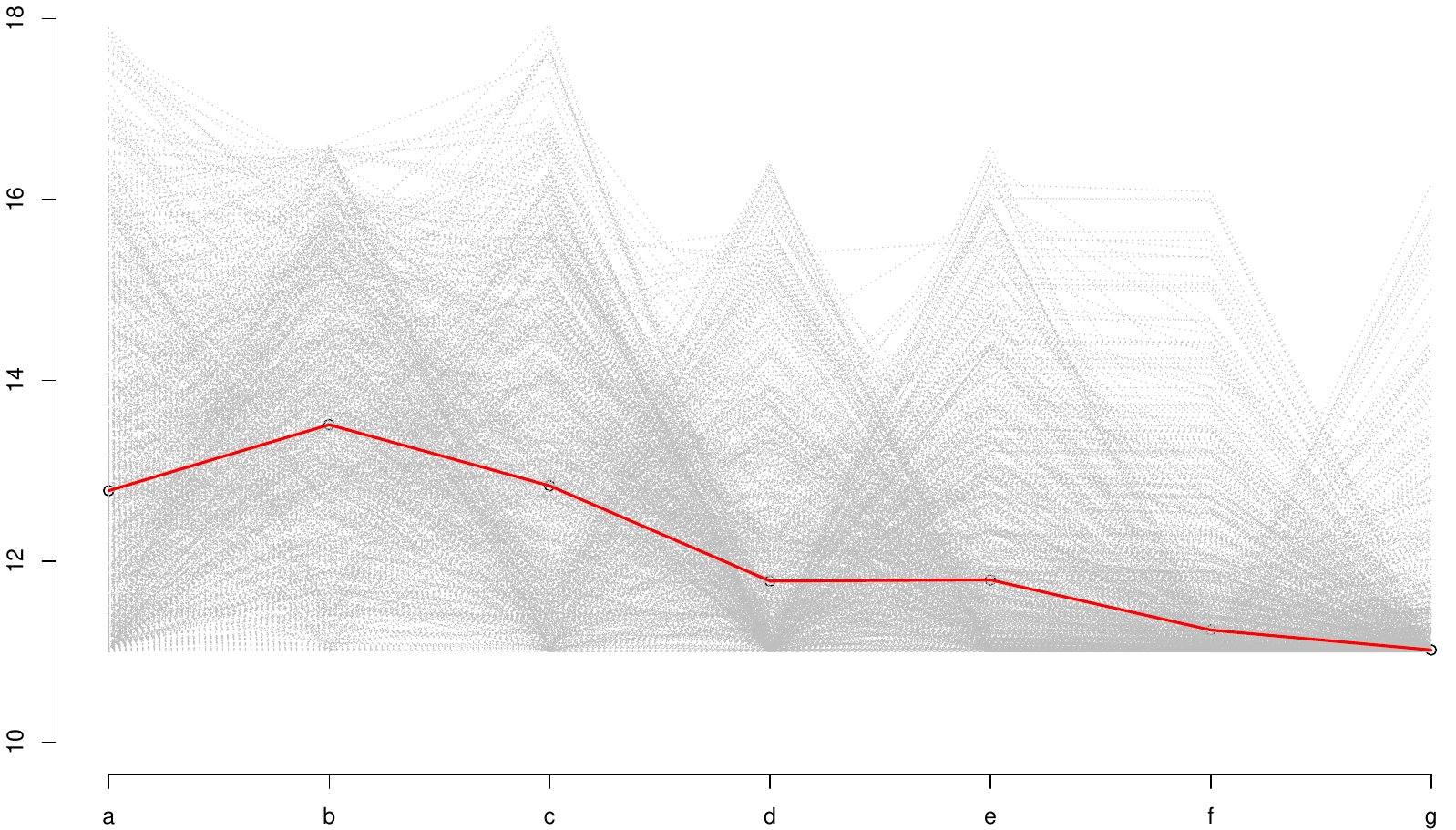}
	\caption{The Choquet integrals of 7 alternatives.} 
	\label{fig-alternatives}
\end{figure}

In fact, graphical tools can also assist us in exploring the nonlinear integral values of various integrands concerning one or more fuzzy measures. For instance, consider the alternatives or integrands presented in Table \ref{table-7-alter-scores}. By utilizing randomly generated 1000 fuzzy measures \cite{wu-NI-quasi-random,wu2020correlative-inoonsistency-random}, we can compute their Choquet integrals and then plot the corresponding lines like in Figure \ref{fig-alternatives}, where the red line represents the median values of their integral values.

\color{black}

\section{Graphic analysis on comparison about fuzzy measures}

\subsection{Comparison analysis of different indices of fuzzy measure}

\begin{figure}[!htp]
	\centering 	\color{black}
	\includegraphics[width=0.6\textwidth]{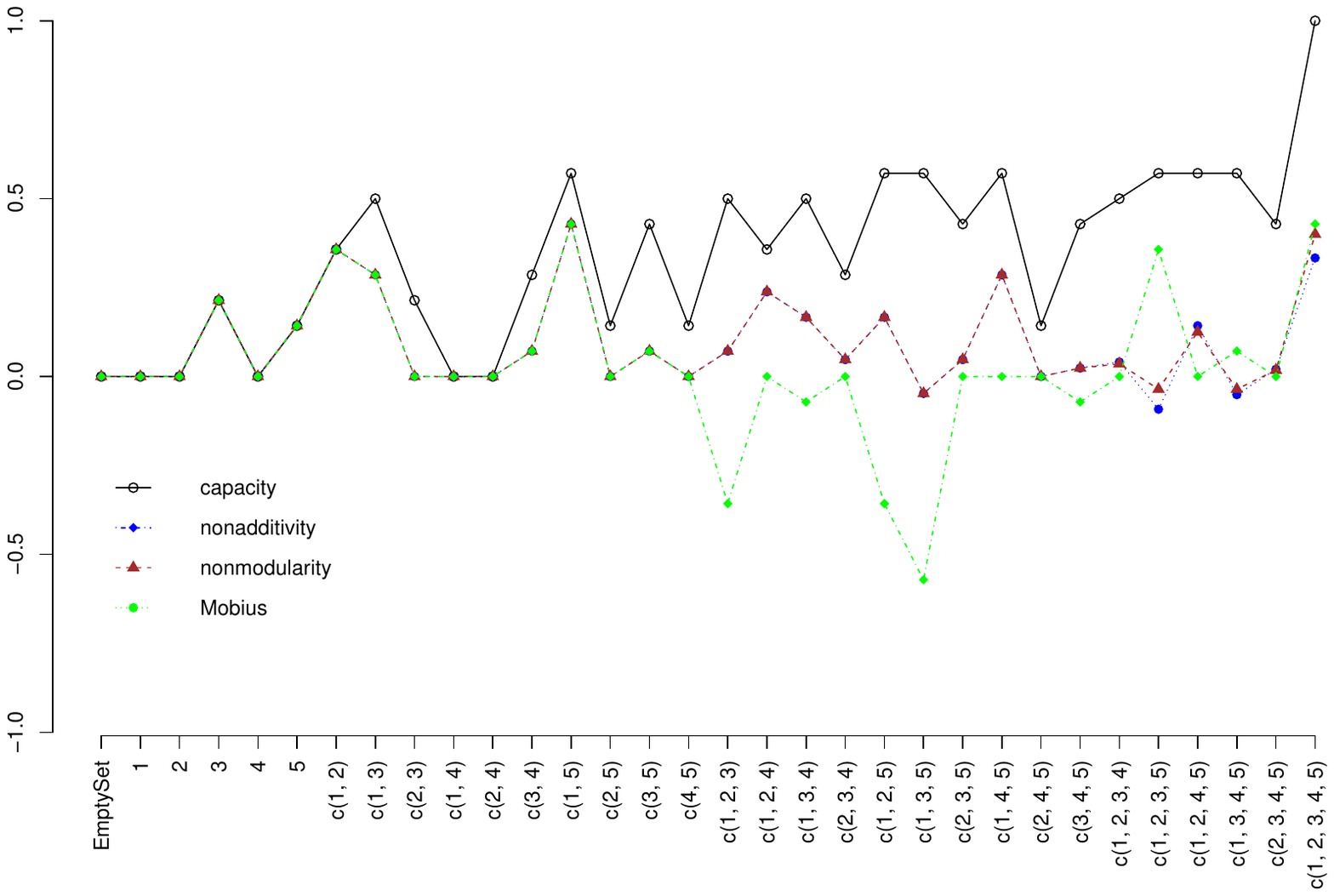}
	\caption{Comparison of three indices for $\mu_{11}$.} 
	\label{fig-2d-of-three-indices}
\end{figure}

For a given fuzzy measure, we can display its related indices in one figure to compare these different indices and discern some characteristics of this fuzzy measure. For example, Figure \ref{fig-2d-of-three-indices} presents three indices of fuzzy measure $\mu_{11}$. It is evident that the nonadditivity and nonmodularity indices are very similar (even matching within 1 $\sim$ 3 orders), while the Möbius representation differs from them.

\subsection{Subsets similarity and clustering analysis}

\begin{figure}[!htp]
	\centering 	\color{black}
	\includegraphics[width=0.6\textwidth]{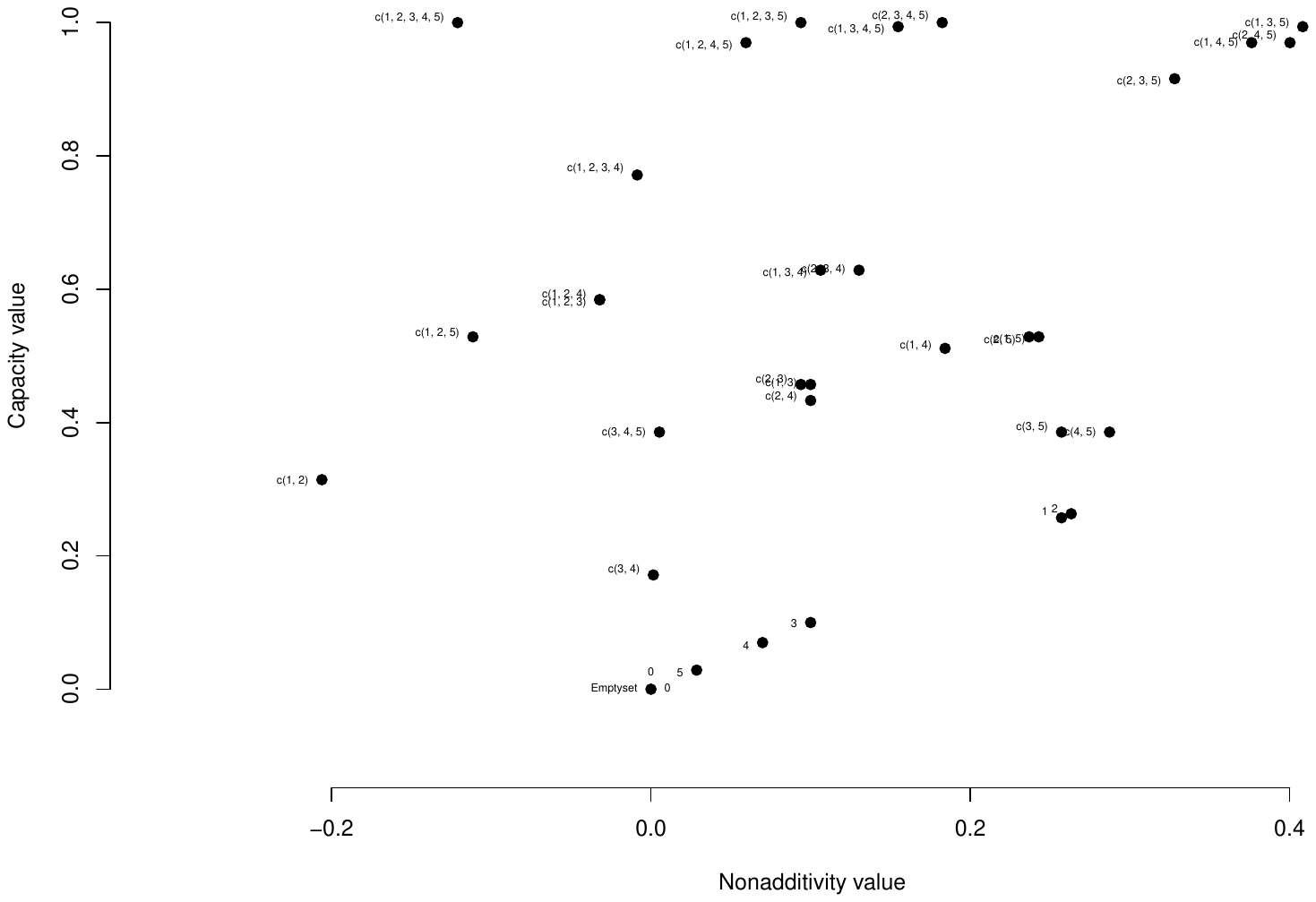}
	\caption{A two dimensional (nonadditivity index, fuzzy measure value) graph of fuzzy measure $\mu_{12}$ for studying subsets.} 
	\label{fig-5-2d-subsets}
\end{figure}

\begin{figure}[!htp]
	\centering 	\color{black}
	\includegraphics[width=1\textwidth]{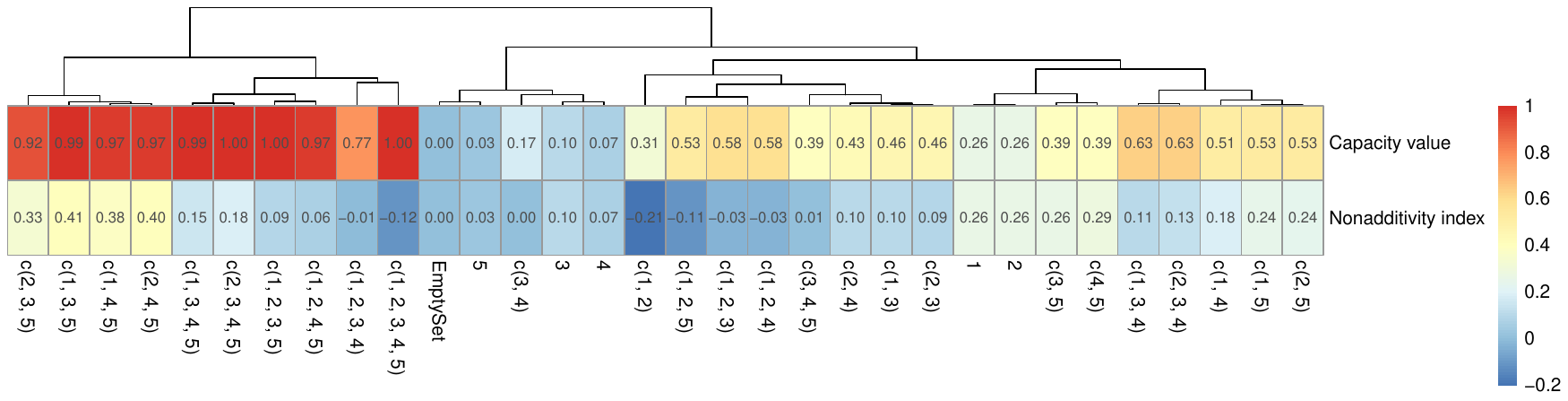}
	\caption{Heatmap for fuzzy measure $\mu_{12}$ of subset clustering with fuzzy measure value and nonadditivity index.} 
	\label{fig-heatmap-5-c-n}
\end{figure}

\begin{figure}[!htp]
	\centering 	\color{black}
	\includegraphics[width=1\textwidth]{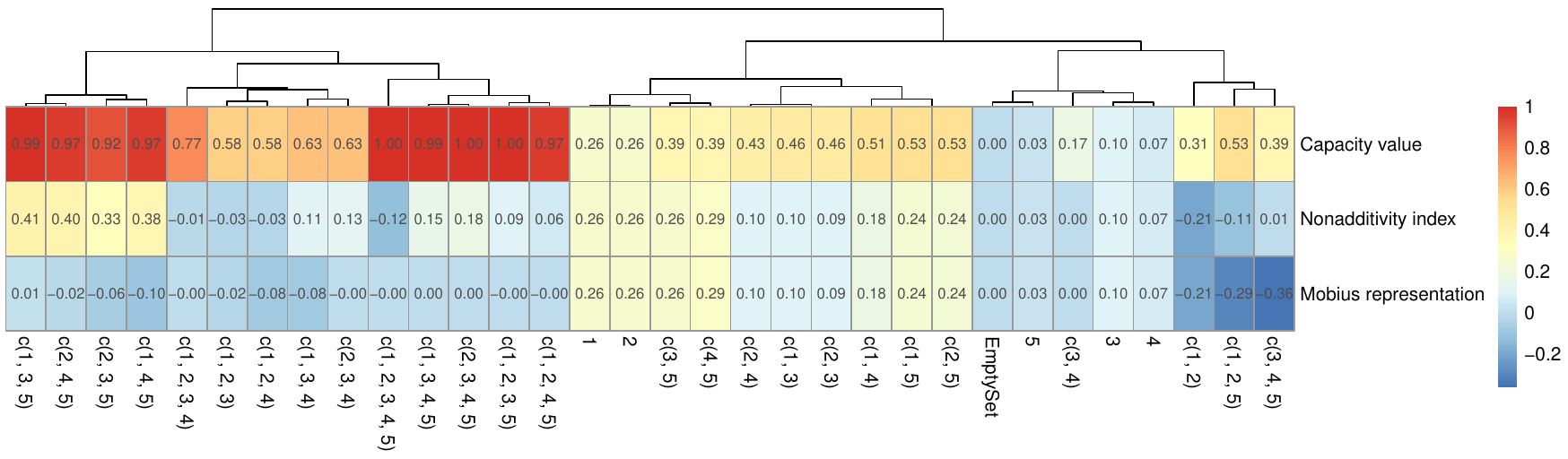}
	\caption{Heatmap for fuzzy measure $\mu_{12}$ of subset clustering further with M\"obius representation.} 
	\label{fig-heatmap-5-c-n-m}
\end{figure}

Using the fuzzy measure as well as its indices, we can analyze the similarity between different subsets and even derive some clusters among them. For instance, for fuzzy measure $\mu_{12}$ on $\{1,2,3,4,5\}$, its fuzzy measure values as well as nonadditivity indices are provided in Figure \ref{fig-heatmap-5-c-n}. We can create a two-dimensional visualization (nonadditivity index, fuzzy measure value) of all the subsets, as shown in Figure \ref{fig-5-2d-subsets}, and cluster them using a heatmap as depicted in Figure \ref{fig-heatmap-5-c-n}. By incorporating the Möbius representation, we obtain a new clustering of these subsets based on three types of indices, presented in Figure \ref{fig-heatmap-5-c-n-m}. It is evident that the heatmap can facilitate clustering analysis for multiple indices of these subsets.



\subsection{Fuzzy measures comparison}
\begin{figure}[!htp]
	\color{black}
	\centering
	\includegraphics[width=0.7\textwidth]{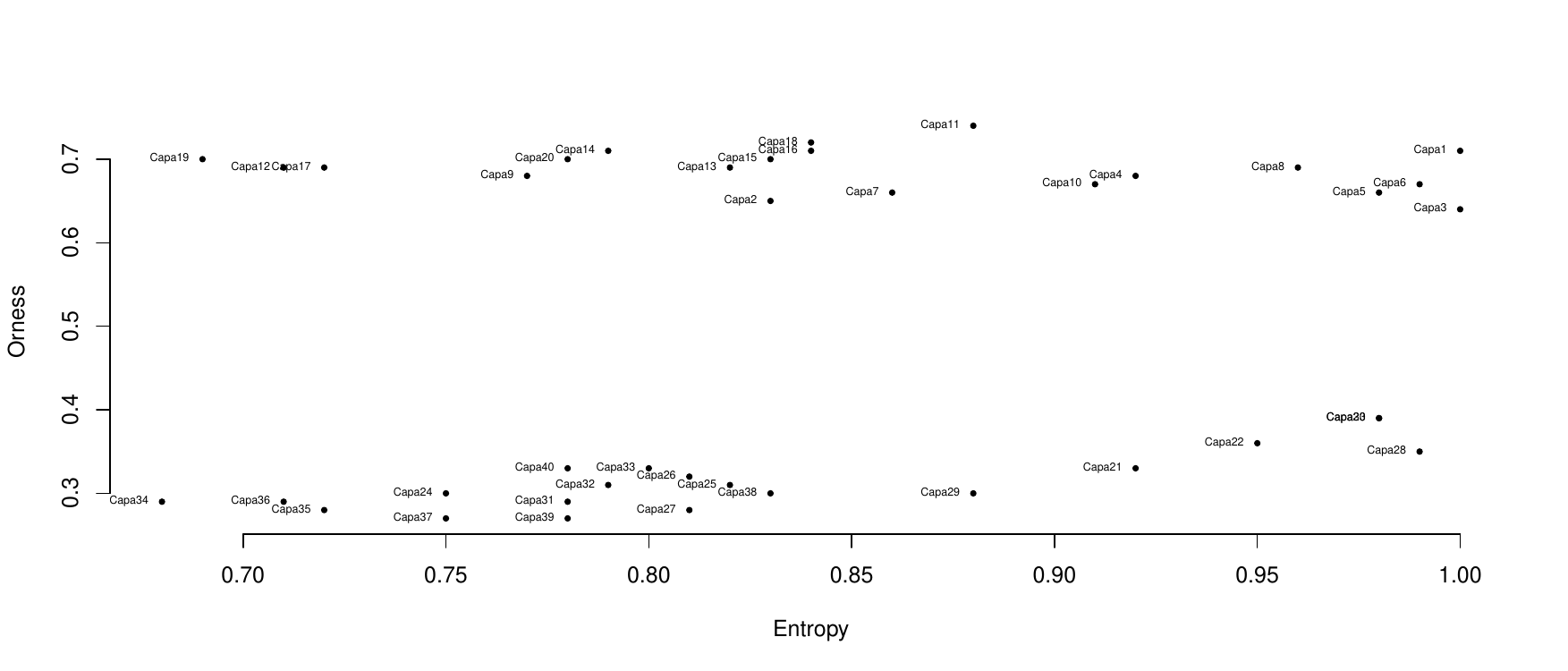}
	\caption{The two dimensional (entropy, ornees) graph of 40 fuzzy measures.} 
	\label{fig-fuzzy measures-2D-40-fuzzy measures}
\end{figure}

\begin{figure}[!htp]
	\color{black}
	\centering
	\includegraphics[width=0.8\textwidth]{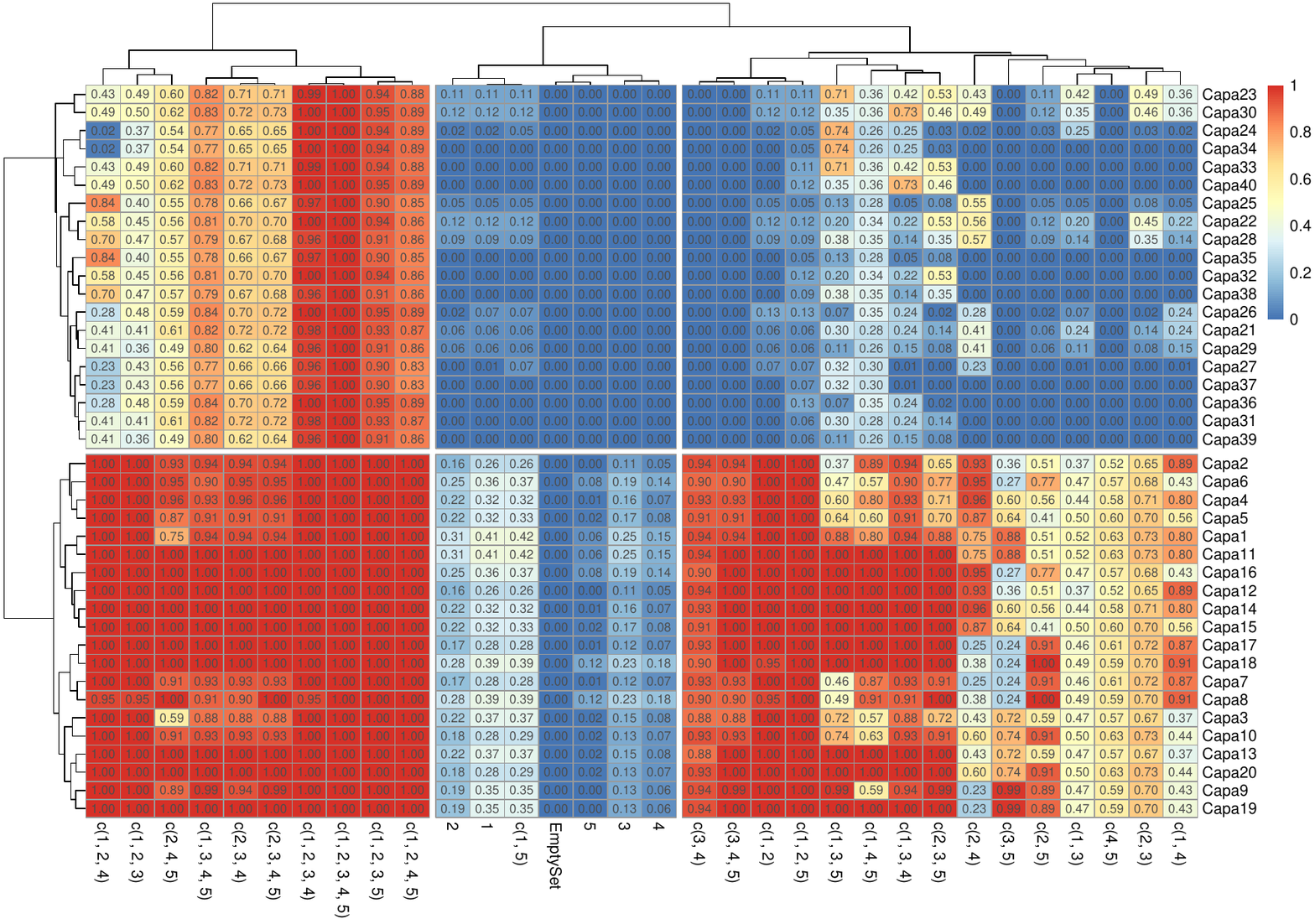}
	\caption{The clustering of two types of special fuzzy measures.} 
	\label{fig-heatmap-for-2-types-40-fuzzy measures}
\end{figure}

Similarly, we can employ the same tools to perform visualization analysis of fuzzy measures. Figure \ref{fig-fuzzy measures-2D-40-fuzzy measures} illustrates the two-dimensional representation of 40 fuzzy measures using entropy and orness indices. Additionally, the clustering of these fuzzy measures based on their fuzzy measure values for all subsets is presented in Figure \ref{fig-heatmap-for-2-types-40-fuzzy measures}.



\section{Conclusions}

Graphic representation provides us with a visualization approach to easily comprehend the concepts, meanings, and relationships within discrete fuzzy measure theory. It has been demonstrated that properties such as monotonicity, duality, additivity, nonadditivity, and nonmodularity can be deduced from the distinctive characteristics or structures of the fuzzy measure graph and can be reflected in the corresponding indices within probabilistic sum terms. Special types of fuzzy measures can be identified through the specific structures or index values within the fuzzy measure. Two-dimensional plots can effectively illustrate the aggregation process of nonlinear integrals, the fuzzy measure fitting procedure for the learning set, as well as comparisons among subsets, fuzzy measures, and integrands. Additionally, heatmap and clustering analyses are valuable for higher-dimensional assessments of subsets and fuzzy measures.

The next major research task may involve the development of interactive and dynamic analysis software or toolkits to facilitate fuzzy measure-based decision-making and machine learning processes.

\bibliographystyle{abbrv}
\bibliography{sample}

\end{document}